\documentclass[a4paper]{amsart}
\usepackage[all]{xy}
\usepackage{times}
\usepackage[bookmarksnumbered,plainpages,hyperindex]{hyperref}
\usepackage{enumerate}
 \newtheorem{theorem}{\sc Theorem}[section]
 \newtheorem{proposition}[theorem]{\sc Proposition}
 \newtheorem{lemma}[theorem]{\sc Lemma}
 \newtheorem{corollary}[theorem]{\sc Corollary}
 \theoremstyle{definition}
 \newtheorem{definition}[theorem]{\sc Definition}
 \newtheorem{definitions}[theorem]{\sc Definitions}
 \newtheorem{example}[theorem]{\sc Example}
 \newtheorem{problem}[theorem]{\sc Problem}
 \theoremstyle{remark}
 \newtheorem{remark}[theorem]{\sc Remark}

 \def\rank{\mathrm{rank}}
 \def\rrank{\mathrm{rrank}}
 \def\Id{\mathrm{Id}}
 \def\Ker{\mathrm{Ker\,}}
 \def\Im{\mathrm{Im\,}}
 \def\Hom{\mathrm{Hom}}
 \def\Alg{\mathrm{Alg}}
 
 \def\End{\mathrm{End}}
 \def\GE{\Gamma_{\boldsymbol{E}}}
 \def\RE{R_{\boldsymbol{E}}}
\def\dtext#1{\emph{\textbf{#1}}}

\begin{document}
\title{On the classification of twisting maps between $K^n$ and $K^m$}


\begin{abstract}
We define the notion of admissible pair for an algebra $A$, consisting on a couple $(\Gamma,R)$, where $\Gamma$ is a quiver and $R$ a unital, splitted and factorizable representation of $\Gamma$, and prove that the set of admissible pairs for $A$ is in one to one correspondence with the points of the variety of twisting maps $\mathcal{T}_A^n:=\mathcal{T}(K^n,A)$. We describe all these representations in the case $A=K^m$.
\end{abstract}

\author{P. Jara}
\thanks{Research partially supported by MTM2007-66666 and FQM-266 (Junta de Andaluc\`{\i}a Research Group). J. L\'{o}pez Pe\~na was supported by Max-Planck Institut f\"ur Mathematik in Bonn and EU Marie-Curie fellowship PIEF-GA-2008-221519. D. \c{S}tefan was financially supported by Contract 560/2009 (CNCSIS code ID$\underline{\phantom{a}}$69).}
\address{Department of Algebra, University of Granada, E-18071--Granada, Spain}
\email{pjara@ugr.es}%

\author{J. L\'opez Pe\~na}
\address{Mathematics Research Centre, Queen Mary University of London, Mile End Road, London E1 4NS, United Kingdom}
\email{jlopez@maths.qmul.ac.uk}%

\author{G. Navarro}
\address{Department of Computer Science and IA, University of Granada, E-18071--Granada, Spain}
\email{gnavarro@ugr.es}%

\author{D. \c Stefan}
\address{University of Bucharest, Faculty of Mathematics, Str. Academiei 14, RO-70109, Bucharest, Romania}
\email{dstefan@al.math.unibuc.ro}%

\date{\today}%
\maketitle

Twisted tensor products of algebras, also known as factorization structures, constitute a particular instance of the notion of \dtext{distributive law} (cf. \cite{Beck69a}) appeared in \cite{Majid90b} and \cite{Tambara90a} and were since studied in different contexts for various purposes, including their applications to braided geometry \cite{Majid90b, Majid91a, Majid94a}, their realization as noncommutative analogues of principal bundles (cf. \cite{Brzezinski98a, Brzezinski00a}), their relation with the (quantum) Yang-Baxter and other nonlinear equations (cf. \cite{VanDaele94a}), and as a suitable replacement for cartesian products in noncommutative geometry (cf. \cite{Cap95a, Jara08a, Lopez08b}).

From a purely algebraic point of view, the notion of twisted tensor product comes directly from the \dtext{factorization problem}:

\begin{quotation} \it
Given some kind of (algebraic) object, is it possible to find two suitable subobjects, having minimal intersection and such that they generate our original object?
\end{quotation}

The factorization problem has been intensively studied in the case of groups, coalgebras, Hopf algebras and algebras (cf. for instance \cite{Agore07a, Caenepeel00a, Caenepeel02a, Takeuchi81a}). In the particular case of algebras, a well known result (independently proven many times) establishes a one-to-one correspondence between the set of factorization structures admitting two given algebras $A$ and $B$ as factors and the set of so-called \dtext{twisting maps}, see Section~\ref{sec:properties}, which allow us to construct the twisted tensor product $A\otimes_{\tau} B$ associated to $\tau$.

Henceforth, the problem of constructing factorization structures with given factors boils down to the problem of finding all the existing twisting maps for those factors. Under suitable, very mild, conditions (for instance, whenever $A$ and $B$ are affine algebras), the set $\mathcal{T}(A,B)$ of all the twisting maps $\tau:B\otimes{A}\to A\otimes{B}$ is an algebraic variety, and two interesting problems arise:

\begin{problem}
Is it possible to describe explicitly the variety $\mathcal{T}(A,B)$?
\end{problem}

\begin{problem}
Once the variety $\mathcal{T}(A,B)$ is known, is it possible to determine which points of the variety give rise to isomorphic algebras?
\end{problem}

These two problems, even in the simplest cases, turn out to be very hard. Though there are many different methods that produce twisted tensor products of two given algebras, not a single one produces all the existing ones is known, let alone describing the properties of the algebraic variety. Even harder is the problem of the determination of the isomorphism classes of algebras obtained from the same factors through different tensor products, or finding any description of these isomorphism classes in terms of the variety $\mathcal{T}(A,B)$.

Recently, Cibils showed in \cite{Cibils06a} that the set $\mathcal{T}(K^{2},A)$ of twisted tensor products between any algebra $A$ and the commutative, semisimple algebra $K^{2}$ (also called the set of \dtext{2--interlaces}) is in one-to-one correspondence with couples of linear endomorphisms of the algebra $A$ satisfying certain conditions. If we take $A=K^n$, these couples of linear maps can be described by combinatorial means using certain families of colored quivers, and this description gives a simple way to describe all the twisted tensor products $K^n\otimes_{\tau} K^{2}$, up to isomorphism (cf. \cite{Cibils06a, Lopez08a}). Some other partial steps in the classification problem for factorization structures have been undertaken in \cite{Borowiec00a} and the final sections of \cite{Guccione99a}.

In the present paper, we extend the combinatorial techniques developed by Cibils, developing the notion of an \dtext{admissible pair} for an algebra $A$, consisting on a couple $(\Gamma, R)$ where $\Gamma$ is a quiver, and $R$ a representation of $\Gamma$ (cf. for instance \cite{Simson, Auslander}) satisfying certain restrictions (namely, to be \dtext{unital}, \dtext{splitted} and \dtext{factorizable}), and prove (cf. Theorem \ref{te:bijection}) that the variety of twisting maps $\mathcal{T}^n_{A}:=\mathcal{T}(K^{n},A)$ is in one to one correspondence with the set of admissible pairs for $A$. Twisted tensor product of the form $A\otimes_{\tau} K^n$ can be reinterpreted as deformations of  the usual tensor product $A \otimes K^n\cong A\times A \times \overset{(n)}{\dotsb} \times A$, which is nothing but the direct product of $n$ copies of the algebra $A$. If $A$ is the algebra of functions defined over the configuration space $Q$ of some mechanical system, then $A\otimes K^n$
 is
  the algebra of functions defined over the configuration space $Q^n$ of the system consisting on $n$ disjoint copies of $Q$. The noncommutative deformations of this algebra obtained as twisted tensor products $A\otimes_{\tau} K^n$ are proposed to serve as toy-model for quantizations of this situation when we assume that the disjoint copies of the physical states are close enough so that they interact with each other.

The paper is structured as follows. In Section \ref{sec:properties}, we study the basic properties of admissible pairs, and introduce the numerical invariants of \dtext{rank} and \dtext{reduced rank} as a measure of the complexity of an admissible pair. We use this notion to characterize the connected components of a quiver in an admissible pair of reduced rank one, proving that splitted, unital and factorizable representations of a quiver $\Gamma$ of reduced rank one are uniquely determined by a set $(M_{i})_{i\in \Gamma^0}$ of two-sided ideals of $A$, and a set $(B_{i})_{i\in \Gamma^0}$ of (unital) subalgebras of $A$ such that

\begin{enumerate}
    \item For each vertex $i\in \Gamma^0$, $A=B_{i}\oplus M_{i}$,
    \item For each arrow $\alpha$ which is not a loop, $M_{s(\alpha)}M_{t(\alpha)}=0$.
\end{enumerate}
A particular example of this setting can be obtained by means of Hochschild extensions of an algebra $B$ with given kernel $M$. Finally, we classify all the splitted, unital, and factorizable representations associated to quivers of reduced rank one, without cycles of length 2, when we take the algebra $A$ to be equal to $K^{m}$.

In Section \ref{sec:absolutelyreducible} we introduce the notion of \dtext{absolutely reducible} (finite dimensional) representations of an algebra $A$, and obtain a canonical form for defining the action of $A$ on an absolutely reducible representation of dimension $n$ by means of a \dtext{normalized invertible matrix}. We pay especial attention to the particular case of two-dimensional absolutely reducible representations, characterizing them all in Theorem \ref{th_rep2}.

Using the aforementioned results, in Section \ref{sec:cyclesoflength2}, we describe all the splitted, unital and factorizable representations of a quiver consisting in a cycle of length two, given a more detailed description for the case of two-dimensional absolutely reducible representations. These results are all merged together in Theorem \ref{th:cycle2}, in which we classify all the representations of a connected quiver of reduced rank one containing a cycle of length 2, thus concluding the classification of all admissible representations of a quiver of reduced rank one.

\section{Twisting maps between $K^n$ and $A$}\label{sec:properties}

Let $A$ be an algebra over a field $K$. The canonical basis of $K^n$ will be denoted by
$\{e_1,\ldots,e_n\}$.  Observe that $\sum_{i=1}^ne_i$ is $1_{K^m}$, the unit of $K^n$. Thus a $K$-linear map
$$
 \tau:K^n\otimes{A}\longrightarrow{A\otimes{K^n}}
$$
is uniquely determined by $\boldsymbol{E}_\tau=(E_{ij})_{i,j=1,\ldots,n}$, a set of
$K$-linear endomorphisms of $A$, satisfying
\begin{equation}\label{eq:E{i,j}}
 \tau(e_i\otimes{a})=\sum_{j=1}^nE_{ij}(a)\otimes{e_j},
 \qquad\forall a\in A,
 \quad\forall{i=1,\ldots,n}.
\end{equation}
Our first aim is to identify the properties that $\boldsymbol{E}_\tau$ has to verify in
order to get a \dtext{twisting map} between the algebras $K^n$ and $A$. Recall that
$\tau:K^n\otimes{A}\longrightarrow{A\otimes{K^n}}$ is a \dtext{twisting map} if, and only if,
the following four conditions hold (cf. \cite{Beck69a, Cap95a}):
\begin{align}
 \label{eq:tau1}
 \tau(m_{K^n}\otimes{A})&=(A\otimes{m_{K^n}})(\tau\otimes{K^n})(K^n\otimes\tau).\\
 \label{eq:tau2}
 \tau(K^n\otimes{m_A})&=(m_A\otimes{K^n})(A\otimes\tau)(\tau\otimes{A}).\\
 \label{eq:tau3}
 \tau(1_{K^n}\otimes{a})&=a\otimes{1_{K^n}},\qquad\forall{a\in{A}}.\\
 \label{eq:tau4}
 \tau(x\otimes{1_A})&=1_A\otimes{x},\qquad\forall{x\in{K^n}}.
\end{align}
In the identities above, $m_{K^n}$ and $m_A$ denote the multiplication on $K^n$ and $A$
respectively, whilst $1_{K^n}$ and $1_A$ denote the units of these algebras.

If we evaluate both sides of relation \eqref{eq:tau1} at $e_i\otimes{e_j}\otimes{a}$, in
view of relation \eqref{eq:E{i,j}}, we get
\begin{align*}
 \tau(e_ie_j\otimes{a})
 &=\sum\limits_{p=1}^n(A\otimes{m_{K^n}})(\tau\otimes{K^n})(e_i\otimes{E_{jp}(a)}\otimes{e_p})\\
 &=\sum\limits_{p=1}^n\sum_{q=1}^nE_{iq}(E_{jp}(a))\otimes{e_pe_q}.
\end{align*}
Therefore
$$
 \sum_{p=1}^n\delta_{i,j}E_{ip}(a)\otimes{e_p}=\sum_{p=1}^n(E_{ip}\circ{E_{jp}})(a)\otimes{e_p}.
$$
In conclusion, relation \eqref{eq:tau1} implies
\begin{equation}\label{eq:tau1p}
 E_{ip}\circ{E_{jp}}=\delta_{i,j}E_{ip},\qquad\forall{i,j,p=1,\ldots,n}.
\end{equation}
Obviously, the converse holds too. Therefore equations \eqref{eq:tau1} and \eqref{eq:tau1p} are equivalent.

The other identities in the definition of twisting maps can be written in a similar way.
By evaluating \eqref{eq:tau2} at $e_i\otimes{a}\otimes{b}$ we get that this relation is
equivalent to
\begin{equation}\label{eq:tau2p}
 E_{ij}(ab)=\sum_{p=1}^nE_{ip}(a)E_{pj}(b),\qquad\forall{i,j=1,\ldots,n},\;\forall{a,b\in{A}}.
\end{equation}
Since the unit of $K^n$ is $\sum_{i=1}^ne_i$, we obtain
\begin{equation}\label{eq:tau3p}
 \sum_{i=1}^nE_{ij}=\Id_A,\qquad\forall{j=1,\ldots,n}.
\end{equation}
Finally, by taking $x:=e_i$ in \eqref{eq:tau4}, we deduce that this relation is equivalent to
\begin{equation}\label{eq:tau4p}
 E_{ij}(1_A)=\delta_{i,j}1_A,\qquad\forall{i,j=1,\ldots,n}.
\end{equation}

In conclusion, if $\mathcal{T}^n_A$ denotes the set of all twisting maps between $K^n$ and
$A$, and $\mathcal{E}^n_A$ denotes the set of systems of endomorphisms
$(E_{ij})_{i,j=1,\ldots,n}$ satisfying relations \eqref{eq:tau1p}--\eqref{eq:tau4p}, we
proved the following theorem.

\begin{proposition}
There is a one-to-one correspondence $\phi_A:\mathcal{T}^n_A\longrightarrow\mathcal{E}^n_A$ that
maps a twisting map $\tau:K^n\otimes{A}\longrightarrow{A}\otimes{K^n}$ to the set of
endomorphisms $\boldsymbol{E}_\tau=(E_{ij})_{i,j=1,\ldots,n}$ which is uniquely defined
such that relation \eqref{eq:E{i,j}} holds true.
\end{proposition}

From now on, instead of working with twisting maps we shall work with systems of
endomorphisms in $\mathcal{E}^n_A$. To such a system
$\boldsymbol{E}=(E_{ij})_{i,j=1,\ldots,n}$ we are going to associate two invariants: a
quiver $\GE$ and a representation $\RE$ of $\GE$.

\begin{definitions}
Let $\Gamma$ be a quiver (cf. \cite{Simson, Auslander}). Let $\Gamma^0$ and $\Gamma^1$ be, respectively, the set of
vertices and the set of arrows of $\Gamma$. The source and the target maps will be
respectively denoted by $s:\Gamma^1\to \Gamma^0$ and $t:\Gamma^1\to \Gamma^0$.
\begin{enumerate}[(1)]
\item
A \emph{path} in $\Gamma$ is a sequence $p=(\alpha_1,\alpha_2,\dots,\alpha_n)$ of arrows such that
$t(\alpha_i)=s(\alpha_{i+1})$, for any $i=1,\dots,n-1$.
\[
\begin{xy}
\xymatrix{
 p:&
 \circ\ar[r]^{\alpha_1}&
 \circ\ar[r]^{\alpha_2}&
 \cdots\ar[r]^{\alpha_{n-1}}&
 \circ\ar[r]^{\alpha_n}&
 \circ
}\end{xy}
\]
The set of paths of length $n$ in $\Gamma$ is denoted by $\Gamma^n$.
\item
For a path $p=(\alpha_1,\alpha_2,\dots,\alpha_n)$ we define the \emph{source} of $p$ by
$s(p):=s(\alpha_1)$. Similarly, the \emph{target} of $p$ is given by $t(p):=t(\alpha_n)$. The  set
of paths $p$ of length $n$ with $s(p)=i$ and $t(p)=j$ is denoted by  $\Gamma^n(i,j)$.
Note that $\Gamma^0(i,i)$ can be identified with the vertex $i$, while $\Gamma^1(i,j)$
is equal to the set of arrows $a$ such that $s(a)=i$ and $t(a)=j$. In particular,
$\Gamma^1(i,i)$ consists of all loops having the vertex $i$ as a source. The source of a
loop will be called a \emph{loop vertex}.
\item
An \emph{oriented cycle} is a path $p$ with $s(p)=t(p)$.
\end{enumerate}
\end{definitions}

\begin{definition}
Let $\boldsymbol{E}:=(E_{ij})_{i,j=1,\dots,n}$ be in $\mathcal{E}^n_A$. The quiver
$\Gamma_{\boldsymbol{E}}$ is defined as follows: the set of vertices of $\GE$ is
$\{v_1,\dots,v_n\}$.  The vertices $v_i$ and $v_j$ are joined by an arrow with the
source in $v_i$ if, and only if, $E_{ij}\neq0$. The vertex $v_i$ shall be represented
simply by $i$ as well.
\end{definition}

For future references we state the following proposition.

\begin{proposition}
The quiver  $\GE$ associated to $\boldsymbol{E}\in\mathcal{E}^n_A$ has no multiple arrows
and, for every vertex $v_i\in\GE^0$, there is a loop having the source (and therefore the
target) in $v_i$.
\end{proposition}
\begin{proof}
By construction, $\GE$ has no multiple arrows. Let $v_i\in\GE^0$. Since $E_{ii}(1_A)=1_A$,
it follows that there is an arrow $\alpha$ such that $s(\alpha)=t(\alpha)=v_i$. This arrow obviously
is unique, as $\GE$ has no multiple arrows.
\end{proof}

Now we are going to construct the second invariant of $\boldsymbol{E}$, namely a certain
linear representation of $\Gamma_{\boldsymbol{E}}$. Recall from \cite{Simson, Auslander} that a
representation of a quiver $\Gamma$ is given by a family of vector spaces
$(V_i)_{i\in\Gamma^0}$ and a family of $K$-linear maps $(\varphi_{\alpha})_{\alpha\in
\Gamma^1}$, where $\varphi_{\alpha}:V_{s(\alpha)}\longrightarrow{V_{t(\alpha)}}$.

\begin{definition}
Let $\Gamma_{\boldsymbol{E}}$ be the quiver associated to
$\boldsymbol{E}=(E_{ij})_{i,j=1,\dots,n}$, a system of endomorphisms in $\mathcal{E}^n_A$.
The representation $R_{\boldsymbol{E}}$ of $\Gamma_{\boldsymbol{E}}$ is defined by the
family of vector spaces $(V_i)_{i\in \Gamma_{\boldsymbol{E}}^0}$ and the family of
$K$-linear maps $(\varphi_{\alpha})_{\alpha\in \GE^1}$, where
$$
 V_i:=A\qquad\text{and}\qquad\varphi_{\alpha}:=E_{s(\alpha),t(\alpha)},
$$
for every $i\in\GE^0$ and $\alpha\in\GE^1$.
\end{definition}

Let $\boldsymbol{E}$ be an element in $\mathcal{E}^n_A$. In the next proposition relations
\eqref{eq:tau1p}--\eqref{eq:tau4p} are rewritten using the maps
$(\varphi_{\alpha})_{\alpha\in\GE^1}$ that define $\RE$.

\begin{proposition}\label{pr:RE}
The maps  $(\varphi_{\alpha})_{\alpha\in \GE^1}$ satisfy the following properties.
\begin{enumerate}[(1)]
\item
For any $i\in \GE^0$, the set $\{\varphi_{\alpha}\mid t(\alpha)=i\}$ is a complete set of
non-zero orthogonal idempotents on $\End_K(A)$. Hence, for any $\alpha',\alpha''\in\GE^1$ such that
$t(\alpha')=t(\alpha'')=i$,
\begin{equation}\tag{$S_i$}\label{split_rep}
 \varphi_{\alpha'}\circ\varphi_{\alpha''}=\delta_{s(\alpha'),s(\alpha'')}\varphi_{\alpha'}
 \qquad
 \text{and}
 \qquad
 \sum_{\{\alpha\mid t(\alpha)=i\}}\varphi_{\alpha}=\Id_A.
\end{equation}
\item
For any arrow $\alpha\in\GE^1$,
\begin{equation}\tag{$U$}\label{unital_rep}
 \varphi_{\alpha}(1_A)=\delta_{s(\alpha),t(\alpha)}1_A.
\end{equation}
\item
For any $i,j\in\GE^0$ and any $a,b\in A$, we have
\begin{equation}\tag{$F_{i,j}$}\label{eq:factorizable}
 \sum_{(\alpha_{1},\alpha_{2})\in \GE^{2}(i,j)}\varphi _{\alpha_{1}}(a)\varphi _{\alpha_{2}}(b)=
 \left\{
 \begin{array}{ll}
 \varphi _{\alpha}(ab), &\text{if}\quad\GE^{1}(i,j)=\{\alpha\}; \\
 0,                 &\text{if}\quad\GE^{1}(i,j)=\emptyset .%
 \end{array}%
 \right.
\end{equation}
\end{enumerate}
\end{proposition}
\begin{proof}
The  relations in $(S_i)$ follow by \eqref{eq:tau1p}, \eqref{eq:tau3p} and the fact that
$E_{kh}=0$, whenever there is no $\alpha\in\GE^1$ such that $s(\alpha)=k$ and $t(\alpha)=h$. In a
similar way $(U)$  follows  from \eqref{eq:tau4p}. Finally, it is not difficult to see
that \eqref{eq:tau2p} is equivalent to the fact that $(F_{i,j})$ holds for every pair
$(i,j)$ of vertices.
\end{proof}

\begin{definition}
Let $A$ be an algebra. We say that $(\Gamma, R)$ is an \emph{admissible pair} of order $n$ for $A$ (shortly an admissible pair) if
\begin{enumerate}[(1)]
\item
$\Gamma$ is a quiver with $n$ vertices such that every vertex is a \emph{loop vertex} and there are no
multiple arrows.
\item
$R$ is a representation of $\Gamma$ such that the vector spaces associated to its
vertices are all equal to $A$. The family $(\varphi_{\alpha})_{\alpha\in \Gamma^1}$ that define $R$
satisfies the following conditions:
\begin{enumerate}[(a)]
\item
The representation $R$ is \emph{splitted}, i.e. relation \eqref{split_rep} holds true
for every $i\in\Gamma^0$.
\item
The representation $R$ is \emph{unital}, that is, relation \eqref{unital_rep} holds
true.
\item
The representation $R$ is \emph{factorizable}, i.e. relation \eqref{eq:factorizable}
holds true for every pair $(i,j)$ of vertices in $\Gamma$.
\end{enumerate}
\end{enumerate}
The set of all admissible pairs $(\Gamma,R)$ of order $n$ for $A$ will be denoted by $\mathcal{R}^n_A$.
\end{definition}

Summarizing, for every $\boldsymbol{E}\in\mathcal{E}^n_A$, we got an admissible pair $(\GE,\RE)$ for $A$ and
$\boldsymbol{E}\mapsto (\GE,\RE)$ defines a map from $\mathcal{E}^n_A$ to $\mathcal{R}^n_A$. By the
proof of Proposition~\ref{pr:RE}, one can see easily that this map is bijective. In
fact, it is sufficient to notice that the inverse maps an admissible pair $(\Gamma,R)$ to
the set $\boldsymbol{E}=(E_{ij})_{i,j=1,\dots,n}$, where $E_{s(\alpha),t(\alpha)}=\varphi_{\alpha}$, for
any arrow $\alpha\in\Gamma^1$, and all other morphisms $E_{ij}$ are zero. In conclusion we have
proved the following.

\begin{theorem}\label{te:bijection}
For an arbitrary $K$-algebra $A$ the sets $\mathcal{T}^n_A$, $\mathcal{E}^n_A$ and
$\mathcal{R}^n_A$ are in one-to-one correspondence.
\end{theorem}

\section{Some basic properties of admissible pairs}

Let $A$ be a $K$-algebra. Throughout this section, $(\Gamma,R)$ will denote an element in
$\mathcal{R}^n_A$. For the family of maps that defines the representation $R$ we shall use the
notation $(\varphi_{\alpha})_{\alpha\in \Gamma^1}$. The maps associated to the loops of $\Gamma$
will play an important r\^{o}le, so we shall use a special notation for them. Namely,
$\varphi_i$ will denote the morphism corresponding to the unique loop $\alpha$ such that
$s(\alpha)=i$.

Our purpose now is to investigate some basic properties of $\Gamma$ and $R$. First, let us
notice that the existence of the representation $R$ imposes some restrictions on $\Gamma$.

\begin{proposition}
Let $\alpha$ be an arrow of $\Gamma$. Then $\varphi_{\alpha}=\Id_A$ if, and only if, $\alpha$ is a loop
and there is no other arrow with target $t(\alpha)$.
\end{proposition}
\begin{proof}
Since $\varphi_{\alpha}$ is the unique non-zero morphism corresponding to an arrow whose target
is $t(\alpha)$, we deduce the required equality by using the second relation in
\eqref{split_rep}. Conversely, let us assume that $\varphi_{\alpha}=\Id_A$. First we prove that
$\alpha$ is the unique arrow having the target in $t(\alpha)$. Let us assume that $\alpha'$ is
another arrow such that $t(\alpha')=t(\alpha)$. Since $\Gamma$ has no multiple arrows, $s(\alpha)\neq
s(\alpha')$. We get
$$
\varphi_{\alpha'}=\varphi_{\alpha'}\circ\varphi_{\alpha}=\delta_{s(\alpha'),s(\alpha)}\varphi_{\alpha}=0.
$$
Since $\varphi_{\alpha'}\neq 0$, by the definition of admissible pair, it follows that $\alpha$ is
the unique arrow having the target in $t(\alpha)$. Again, by the definition of admissible
pair, every vertex is a loop vertex. In the view of the foregoing, $\alpha$ has to be a
loop.
\end{proof}

By the above Proposition, for every vertex of $\Gamma$, there are two possibilities, as it
is indicated in the two pictures below.
$$
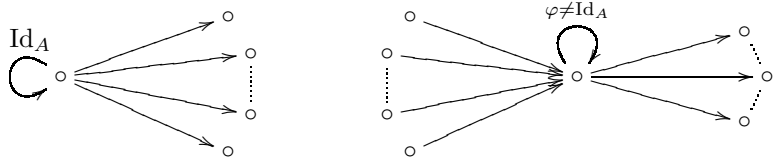
\xy
 (0,10)*+{{\circ}}="a",(-4,15)*+{\Id_A},
 (22,18)*+{{\circ}}="b",
 (25,13)*+{{\circ}}="c",(25,5)*+{{\circ}}="c1",
 (22,0)*+{{\circ}}="e",
 \ar@(ul,dl) "a";"a"
 \ar "a";"b"
 \ar "a";"c"\ar "a";"c1"
 \ar "a";"e" \ar@{..} "c";"c1"
\endxy
\qquad\qquad
 \xy
 (0,10)*+{{\circ}}="a",
 (22,16)*+{{\circ}}="b",
 (25,10)*+{{\circ}}="c",
 (22,4)*+{{\circ}}="e",
 (-22,18)*+{{\circ}}="f",
 (-25,13)*+{{\circ}}="g",
 (-22,0)*+{{\circ}}="k",(-25,5)*+{{\circ}}="k1",
 \ar "a";"b"
 \ar "a";"c"
 \ar "a";"e"
 \ar "f";"a"
 \ar "g";"a"
 \ar"k";"a"\ar"k1";"a" \ar@{..}"k1";"g"\ar@{..}"e";"c"\ar@{..}"b";"c"
 \ar@(ul,ur) "a";"a"^{\varphi\neq\Id_A}
\endxy
$$
In the first one it is represented the case when the morphism corresponding to the loop
is the identity of $A$. In the second picture, the morphism $\varphi$,
associated to the loop, is not the identity of $A$. In what follows in the pictures we
only draw the loops $\alpha$ with $\varphi=\mathrm{Id}_A$.

\vspace*{1ex}

In order to measure the complexity of $\Gamma$ we introduce a numerical invariant, the \dtext{rank}
of a vertex. Let $i\in\Gamma^0$. We set
\begin{equation}\label{eq:rank}
 \rank(i):=\sum_{j}\#\Gamma^1(j,i)=\#\{\alpha\in\Gamma^1\mid t(\alpha)=i\}.
\end{equation}
We also define the \dtext{reduced rank} of $i$ as
\begin{equation}\label{ep:rrank}
 \rrank(i):=\sum_{j\neq i}\#\Gamma^1(j,i)=\#\{\alpha\in\Gamma^1\mid t(\alpha)=i\text{ and }s(\alpha)\neq i\}.
\end{equation}

Note that, in the quiver of an admisible pair,  $\rrank(i)=0$ if, and only if, there is only one arrow $\alpha$ with $t(\alpha)=i$. Of course, in this case we also have $s(\alpha)=i$ and $\varphi_{\alpha}=\Id_A$. If $\rrank(i)=r$, then there are exactly $r$ arrows having their target in $i$ and which are not loops.

\begin{definition}
The \dtext{rank} of $\Gamma$ is defined by
\begin{equation}\label{eq:rankg}
 \rank(\Gamma):=\max\{\rank(i)\mid\;i\in{\Gamma^0}\}.
\end{equation}
In the same way, we may define the \dtext{reduced rank} of $\Gamma$.
\end{definition}

\begin{proposition}\label{pr:rank}
Let $A$ be a $K$-algebra of dimension $m$. If $\boldsymbol{E}\in\mathcal{E}^n_A$, then $\rrank(\Gamma_{\boldsymbol{E}})\leq \mathrm{min}(n-1,m-1)$.
\end{proposition}
\begin{proof}
Obviously, $\rrank(\Gamma_{\boldsymbol{E}})\leq n-1$, as there are at most $n-1$ arrows $\alpha$ such
that $t(\alpha)=i$ and $s(i)\neq{i}$, for any vertex $i$ in $\Gamma_{\boldsymbol{E}}$. On the other
hand, for $i\in\Gamma^0_{\boldsymbol{E}}$, let $X=\{\alpha\in \Gamma^1_{\boldsymbol{E}}\mid t(\alpha)=i\}$, then
$\{\varphi_{\alpha}\mid\alpha\in{X}\}$ is a complete set of orthogonal idempotents. Therefore $A$
decomposes as a direct sum of non zero vector subspaces $A=\bigoplus_{\alpha\in X}W_{\alpha}$,
where $W_{\alpha}=\Im\varphi_{\alpha}$. Hence
$$
\# X\leq\displaystyle\sum_{\alpha\in X}\mathrm{dim}W_{\alpha}=\mathrm{dim}A.
$$
It results $\rrank(i)\leq\mathrm{dim} A-1$. Thus the proposition is proved.
\end{proof}

The simplest non-trivial quivers (i.e. containing arrows that are not loops) are those
of reduced rank $1$. Their connected components are described in the following
proposition. Note that when we speak about the connected components of $\Gamma$, we mean the
connected components of the undirected graph obtained by removing the orientation of all
arrows in $\Gamma$. On the other hand, each connected component can be seen as a quiver with
respect to the orientation of its edges that is inherited from $\Gamma$.

\begin{proposition}\label{pr:24}
Let $\Gamma$ be an \emph{arbitrary} finite quiver of reduced rank one, then any connected component contains at most a unique cycle which is not a loop.
\end{proposition}
\begin{proof}
Clearly we may assume that the connected component $\Lambda$ is not reduced to a unique vertex.
Let us show that in the same component cannot exist two different cycles. We pick up a vertex $i$
in the first cycle and a vertex $j\neq{i}$ in the second one.
Thus, as $\Lambda$ is connected, there is a sequence of vertices $i=i_1,\dots,i_{h+1}=j$ and a
sequence of arrows $\alpha_1,\dots,\alpha_{h}$ in $\Lambda$ such that either $s(\alpha_k)=i_k$ and
$t(\alpha_k)=i_{k+1}$ or $s(\alpha_k)=i_{k+1}$ and $t(\alpha_k)=i_k$, for any $k=1,\dots,h$.
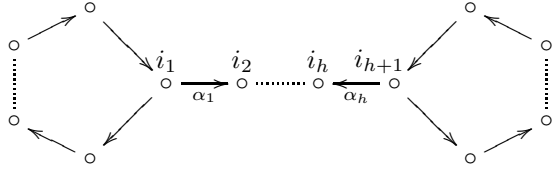
$$
\xy
 (0,-10)*+{\circ}="a",(0,10)*+{\circ}="aa",
 (-10,-5)*+{\circ}="a1",(-10,5)*+{\circ}="aa1",
 (60,-5)*+{\circ}="b1",(60,5)*+{\circ}="bb1",
 (50,-10)*+{\circ}="b",(50,10)*+{\circ}="bb",
 (10,0)*+{\circ}="i",(10,3)*+{i_1}="i1",
 (20,0)*+{\circ}="k",(20,3)*+{i_2}="k1",
 (30,0)*+{\circ}="l",(30,3)*+{i_h}="l1",
 (40,0)*+{\circ}="j",(38,3)*+{i_{h+1}}="j1",
 \ar "aa";"i" \ar "i";"a" \ar@{..} "a1";"aa1"\ar@{..} "b1";"bb1"
 \ar@{->} "aa1";"aa" \ar@{->} "a";"a1"\ar@{->} "b";"b1"\ar@{->} "bb1";"bb"
 \ar@{->} "i";"k"_{\alpha_1}
 \ar@{..} "k";"l"
 \ar@{->} "j";"l"^{\alpha_h}
 \ar "j";"b" \ar "bb";"j"
\endxy
$$
Since $\rrank(i)=1$ it follows that $s(\alpha_1)=i_1=i$ and $t(\alpha_1)=i_2$. By induction, as
$\rrank(i_k)=1$, we get $s(\alpha_k)=i_k$ and $t(\alpha_{k})=i_{k+1}$, for any $k=1,\dots,n$. It
follows that $t(\alpha_{h})=i_{h+1}=j$, so $\rrank(i_{h+1})\geq 2$, which is impossible.

\end{proof}

\begin{proposition}
With the same assumption as in Proposition~\ref{pr:24}, every vertex in the possible oriented cycle is the root of a (naturally) oriented tree which can be degenerated (i.e., with only one vertex) and such that two different of these trees are disjoint.
\end{proposition}
\begin{proof}
Let us denote by $\Lambda'$ the quiver obtained removing all the arrows $\alpha_1,\ldots,\alpha_r$ in the cycle with vertices $i_1,\ldots,i_r$ such that $s(\alpha_k)=i_k$ for any $k=1,\ldots,r$. See the figure below.
$$
 
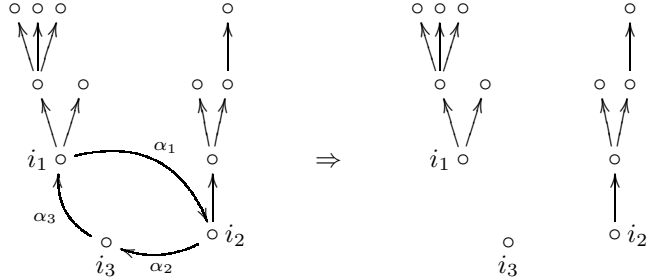
\xy
 (0,0)*+{\circ}="a1",(6,-11)*+{\circ}="a3",(20,-10)*+{\circ}="a2",
 (-3,0)*+{i_1},(6,-14)*+{i_3},(23,-10)*+{i_2},
 (-3,10)*+{\circ}="b",
 (3,10)*+{\circ}="c",
 (-6,20)*+{\circ}="d",
 (-3,20)*+{\circ}="e",
 (0,20)*+{\circ}="f",
 (20,0)*+{\circ}="h",
 (18,10)*+{\circ}="k",(22,10)*+{\circ}="l",(22,20)*+{\circ}="m",
 \ar@/^3.5ex/ "a1";"a2"^{\alpha_1}\ar@/^1.3ex/"a2";"a3"^{\alpha_2}\ar@/^1.7ex/
 "a3";"a1"^{\alpha_3}
 \ar "a1";"b"\ar "a1";"c"\ar "b";"d"\ar "b";"e"\ar "b";"f"
 \ar "a2";"h"
 \ar "h";"k" \ar "h";"l" \ar "l";"m"
 \endxy
\qquad\Rightarrow\qquad
 \xy
 (0,0)*+{\circ}="a1",(6,-11)*+{\circ}="a3",(20,-10)*+{\circ}="a2",
 (-3,10)*+{\circ}="b",(-3,0)*+{i_1},(6,-14)*+{i_3},(23,-10)*+{i_2},
 (3,10)*+{\circ}="c",
 (-6,20)*+{\circ}="d",
 (-3,20)*+{\circ}="e",
 (0,20)*+{\circ}="f",
 (20,0)*+{\circ}="h",
 (18,10)*+{\circ}="k",(22,10)*+{\circ}="l",(22,20)*+{\circ}="m",
 \ar "a1";"b"\ar "a1";"c"\ar "b";"d"\ar "b";"e"\ar "b";"f"
 \ar "a2";"h"
 \ar "h";"k" \ar "h";"l" \ar "l";"m"
 \endxy
$$
We fix $k\in\{1,\dots,r\}$. The connected component $\Lambda'_k\subseteq\Lambda'$ that contains $i_k$ has no cycles, otherwise there would be two different cycles in $\Lambda$. Hence, $\Lambda'_k$ is a tree (recall that a graph is a tree if, and only if, it is connected and does not contains cycles, or equivalently,
any two vertices are connected by a unique path). Since the rank of $\Lambda$ is $1$, the arrows of $\Lambda'_k$ are oriented in the canonical way, that is, the one such that the arrows at the root are outgoing.
\end{proof}

\begin{lemma}\label{le:factorizable}
Let $(\Gamma,R)$ be an admissible pair.
\begin{enumerate}[(1)]
\item
If $i\in{\Gamma^0}$ is not a vertex in a cycle of length $2$ then condition $(F_{i,i})$ is equivalent to the fact that $\varphi_i$ is an algebra map (recall that $\varphi_i$ denotes the morphism corresponding to the loop having the source in $i$).
\item
Let $(\alpha_1,\alpha_2)\in\Gamma^2$ be a path such that $\alpha_2$ and $\alpha_1$ are as in the following picture,
$$
\xy
 (0,0)*+{{\circ}}="i",(0,3)*+{i},
 (26,0)*+{{\circ}}="j",(26,3)*+{k},
 (13,0)*+{{\circ}}="k",(13,3)*+{j},
 \ar"i";"k"^{\alpha_1} \ar"k";"j"^{\alpha_2}
\endxy
$$
If $\rrank(j)=\rrank(k)=1$, then conditions $(F_{i,k})$ and $(F_{j,k})$ are equivalent. Moreover, these conditions are also equivalent to
\begin{equation}\label{eq:Fjk}
 \Ker(\varphi_{s(\alpha_2)})\Ker(\varphi_{t(\alpha_2)})=0.
\end{equation}
\end{enumerate}
\end{lemma}
\begin{proof}
By assumption $i$ is not a vertex in a cycle of length $2$. Then, either $\Gamma^1(i,k)=\emptyset$ or $\Gamma^1(k,i)=\emptyset$, for any $k\neq i$. Therefore, in this case, equation $(F_{i,i})$ becomes
\[
\varphi_i(a)\varphi_i(b)=\varphi_i(ab),\qquad\forall{a,b\in{A}}.
\]
Let us prove the second claim. We first show that $k$ is not a vertex in a cycle of length 2. Indeed, if $l$ were the second vertex of such a cycle and $l\neq j$ then we would have $\rrank(k)\geq 2$. Thus $l=j$. On the other hand, this equality implies $\rrank(j)>1$, which is also impossible.
Since $k$ is not the vertex of a cycle of length 2, it follows that $\varphi_k$ is an algebra map.
Since $\Gamma^2(i,k)=\{(\alpha_1,\alpha_2)\}$ and there is no arrow having the source in $i$ and the target in $k$, relation $(F_{i,k})$ is equivalent to
\[
 \varphi_{\alpha_1}(a)\varphi_{\alpha_2}(b)=0,\qquad\forall{a,b\in{A}}.
\]
By hypothesis, the reduced rank of $j$ and of $k$ is one. Therefore, $\varphi_{\alpha_1}=\Id_A-\varphi_j$ and  $\varphi_{\alpha_2}=\Id_A-\varphi_k.$ Then the above equality is equivalent to
\begin{equation}\label{eq:Fjk'}
 [a-\varphi_j(a)][b-\varphi_k(b)]=0,\qquad\forall{a,b\in{A}}.
\end{equation}
Obviously, $\Gamma^2(j,k)=\{(\lambda_1,\alpha_2),(\alpha_2,\lambda_2)\}$, where $\lambda_1$ and $\lambda_2$ are the unique loops such that $s(\lambda_1)=j$ and $s(\lambda_2)=k$. As $\varphi_{\lambda_1}=\varphi_j$ and  $\varphi_{\lambda_2}=\varphi_k$, relation $(F_{j,k})$ can be written as follows
\begin{equation}\label{eq:Fjk''}
 \varphi_j(a)\varphi_{\alpha_2}(b)+\varphi_{\alpha_2}(a)\varphi_k(b)=\varphi_{\alpha_2}(ab),\qquad\forall{x,y\in{A}}.
\end{equation}
Since $\varphi_{\alpha_2}=\Id_A-\varphi_k$ one can prove easily that \eqref{eq:Fjk'} and \eqref{eq:Fjk''} are equivalent. In conclusion, $(F_{i,k})$ and$(F_{j,k})$ are equivalent too. To conclude it is enough to prove that both are equivalent to \eqref{eq:Fjk}. Relation \eqref{eq:Fjk'} can be written as:
\[
\Im(\Id_A-\varphi_{s(\alpha_2)})\Im(\Id_A-\varphi_{t(\alpha_2)})=0.
 \]
On the other hand, $\Im(\Id_A-\varphi_{s(\alpha_2)})=\Ker(\varphi_{s(\alpha_2)})$, as $\varphi_{s(\alpha_2)}$ is an idempotent $K$-linear map. A similar equality holds for $\varphi_{t(\alpha_2)}$, so the lemma is proved.
\end{proof}

\begin{theorem}\label{te:Rank1Gen}
Let $A$ be a $K$-algebra and let $\Gamma$ be a quiver such that its vertices are loop vertices and it has no multiple arrows. Assume that $\rrank(\Gamma)=1$ and that $\Gamma$ does not contain any cycle of length $2$. Then to give a splitted, unital and factorizable representation of $\Gamma$ over $A$ is equivalent to give a set of idempotent algebra
endomorphisms $(\varphi_i)_{i\in\Gamma^0}$ such that, for any arrow $\alpha:i\to{j}$ which is not a loop,
\begin{equation}\label{eq:phi(ab)}
 \Ker\varphi_{i}\Ker\varphi_{j}=0.
\end{equation}
\end{theorem}
\begin{proof}
Let $R$ be a splitted, unital and factorizable representation of $\Gamma$. Let $(\varphi_{\alpha})_{\alpha\in\Gamma^1}$ denote the family of $K$-linear endomorphisms that defines $R$. For $i\in\Gamma^0$ we take $\varphi_i$ to be the morphism corresponding to the loop $\alpha$ that has the source in $i$. Since $\Gamma$ does not contain cycles of length $2$, by the previous lemma, it follows that $\varphi_{i}$ is an algebra map, for any $i\in{\Gamma}^0$. Let $\alpha\in\Gamma^1$ which is not a loop. If $\varphi_{s(\alpha)}=\Id_A$ we have nothing to prove.

Let us consider the case when $\varphi_{s(\alpha)}\neq\Id_A$. Thus there is an arrow $\beta$, which is not a loop, such that $t(\beta)=s(\alpha)$. Hence relation \eqref{eq:phi(ab)} follows by the second part of Lemma~\ref{le:factorizable}.

Conversely, let  $(\varphi_i)_{i\in\Gamma^0}$ be a family of idempotent algebra morphisms satisfying relation
\eqref{eq:phi(ab)}. We want to construct a splitted, unital and factorizable representation $R$ of $\Gamma$. For every arrow $\alpha$ which is not a loop, we set $\varphi_{\alpha}=\Id_A-\varphi_{t(\alpha)}$. Obviously, $\{\varphi_{\alpha},\varphi_{t(\alpha)}\}$ is a complete set of orthogonal idempotents, so the representation $R$ defined by $(\varphi_{\alpha})_{\alpha\in\Gamma^1}$ is splitted. Trivially $R$ is unital, as $\varphi_i$ is a morphism of algebras, for any $i\in\Gamma^0$. It remains to prove that $R$ is factorizable. Since $\rrank(\Gamma)=1$ and $\Gamma$ has not cycles of length 2, the non-trivial relations that can occur are $(F_{i,k})$ and $(F_{j,k})$, where the vertices $i$, $j$ and $k$ are as in the picture below.
$$
\xy
 (0,0)*+{{\circ}}="i",(0,3)*+{i},
 (26,0)*+{{\circ}}="j",(26,3)*+{k},
 (13,0)*+{{\circ}}="k",(13,3)*+{j},
 \ar"i";"k"^{\beta} \ar"k";"j"^{\alpha}
\endxy
$$
By Lemma~\ref{le:factorizable} (2) these relations are equivalent to \eqref{eq:phi(ab)}, so the Theorem is proved.
\end{proof}

\begin{corollary}
Let $A$ be a finite dimensional $K$-algebra and let $\Gamma$ be a quiver of reduced rank one without cycles of length $2$. A splitted, unital and factorizable representation of $\Gamma$ is uniquely defined by a set
$(M_i)_{i\in\Gamma^0}$ of two-sided ideals and a set $(B_i)_{i\in\Gamma^0}$ of unital subalgebras of $A$ that satisfy the following two conditions:
\begin{enumerate}[(1)]
\item
For each vertex $i\in\Gamma^0$, $A=B_i\bigoplus M_i$;
\item
For each arrow $\alpha$ which is not a loop, $M_{s(\alpha)}M_{t(\alpha)}=0$.
\end{enumerate}
\end{corollary}
\begin{proof}
For a vertex $i$, take in the previous theorem $B_i:=\Im\varphi_i$ and $M_i:=\Ker\varphi_i$. Conversely, for a family of ideals $(M_i)_{i\in\Gamma^0}$, satisfying the conditions in the corollary, we define $\varphi_i:=\sigma_i\circ\pi_i$, where $\sigma_i$ is the inclusion $B_i\subseteq A$ and $\pi_i$ is the projection of $A=B_i\bigoplus M_i$ onto $B_i$.
\end{proof}

\begin{remark} \label{re:rank1}
Let $A$ be a $K$-algebra. We assume that there is a two-sided ideal $M\leq A$ such that $M^2=0$. Let
$\Gamma$ be a quiver of reduced rank one and let $(\varphi_i)_{i\in\Gamma^0}$ be a family of idempotent algebra endomorphisms of $A$. If
$$
 \Ker\varphi_i\subseteq M
$$
for any $i\in\Gamma^0$, then condition \eqref{eq:phi(ab)} is automatically satisfied. Hence the family $(\varphi_i)_{i\in\Gamma^0}$ induces a splitted, unital and factorizable representation of $\Gamma$.
\end{remark}

In order to give examples of families $(\varphi_i)_{i\in \Gamma^0}$ satisfying the assumptions in the above remark, we recall the construction of \dtext{Hochschild extensions} (also called in the literature \dtext{abelian extension}, cf. \cite[Sec. 1.5.3]{Loday98a}) of an algebra $B$ with kernel a given $B$-bimodule $M$. By definition, such an extension is given by a normalized Hochschild 2-cocycle, that is, a
$K$-linear map $\omega:B\otimes B\rightarrow M$ that verifies the following equalities:
\begin{itemize}
\item
$\omega(a\otimes1_B)=\omega(1_B\otimes a)=0$, for all $a\in B$.
\item
$a\cdot\omega(b\otimes c)-\omega(ab\otimes c)+\omega(a\otimes{bc})-\omega(a\otimes b)\cdot c=0$, for all $a,b,c\in{B}$.
\end{itemize}

To these data one can associate a unital and associative algebra $A$ as follows. As a vector space $A:=B\oplus{M}$. The multiplication on $A$ is defined by
$$
 (b,m)\cdot (b',m')=(bb',bm'+mb'+\omega(b\otimes b'))
$$
and the unit is $(1,0)$. Note that $M$ is a two-sided ideal in $A$ and $M^2=0$ (of course $M$ can be identified with a subset of $A$ via the map $m\mapsto (0,m)$).

\begin{lemma}
Let $f:M\rightarrow M$ be a $K$-linear map.
\begin{enumerate}[(1)]
\item
The $K$-linear map $\varphi:A \rightarrow A$, defined by $\varphi(b,m)=(b,f(m))$ for any $b\in B$ and $m\in M$, is an algebra map if, and only if, $f$ is a morphism of $B$-bimodules and $f\circ\omega=\omega$.
\item
$\varphi$ is idempotent if, and only if, $f$ is so.
\item
The kernel of $\varphi$ is contained into $M$.
\end{enumerate}
\end{lemma}
\begin{proof}
$\varphi$ is a morphism of algebras if, and only, if
\begin{equation}\label{eq:f}
 f(bm'+mb')+(f\circ \omega)(b\otimes b')=b f(m')+f(m)b'+\omega(b\otimes b')
\end{equation}
for any $b,b'\in B$ and $m,m'\in M$.

Let us assume that $\varphi$ is a morphism of algebras, hence the above equation holds. By taking $b'=0$ we get that $f$ is a morphism of left $B$-modules. Similarly, we deduce that $f$ is a map of right $B$-modules. In particular,
$$
f(mb'+bm')=f(m)b'+bf(m'),
$$
so $f\circ\omega=\omega$ holds too. Conversely, if $f$ is a morphism of $B$-bimodules and $f\circ\omega=\omega$, then $f$ obviously satisfies the equation \eqref{eq:f}, so $\varphi$ is a morphism of algebras.

The second and the third part of the Lemma are obvious.
\end{proof}

\begin{theorem}
Let $A$ be the Hochschild extension associated to $(B,M,\omega)$, where $B$ is a $K$-algebra, $M$ is a $B$-bimodule and $\omega$ is a 2-cocycle of $B$ with coefficients in $M$. Let $\Gamma$ be a quiver of reduced rank one without cycles of length two. If $(f_i)_{i\in\Gamma^0}$ is a family of idempotent endomorphisms of $M$ and $\Im\,\omega\subseteq \bigcap_{i\in \Gamma^0} \Im f_i$ then $(\varphi)_{i\in \Gamma^0}$ induces a splitted, unital and factorizable representation of $\Gamma$, where $\varphi_i$ is constructed from $f_i$ as in the previous Lemma.
\end{theorem}
\begin{proof}
It is a direct consequence of Theorem~\ref{te:Rank1Gen} and the previous Lemma. See also Remark~\ref{re:rank1}.
\end{proof}

Throughout the remaining of this section we take $A=K^m$. Our purpose is to classify all splitted, unital and factorizable representations (over $A$) of a quiver $\Gamma$ of reduced rank one and without cycles of length two.

Let $\theta:K^m \rightarrow K^m$ be an algebra endomorphism of $K^m$. If $\{f_1,\ldots, f_m\}$ is the canonical basis on $K^m$, then $\theta(f_p)$ is an idempotent of $K^m$. Therefore there is a set $\Theta_p\subseteq \{1,\ldots , m\}$ such that
$$
\theta(f_p)=\displaystyle \sum_{q\in \Theta_p} f_q.
$$
Note that if $\Theta_p=\emptyset$ then $\theta(f_p)=0$. Thus the kernel of $\theta$ is the vector subspace of $K^m$ generated by all elements $f_p$ such that $\Theta_p=\emptyset$. Moreover, $(\Theta_p)_{p=1,\ldots ,m}$ is a partition of $\{1,\ldots, m\}$ in $m$ (possibly empty) subsets. Indeed, if $p\neq q$ then $f_pf_q=0$, so $$
 0=\theta(f_pf_q)
 =\sum_{r\in \Theta_p}\sum_{s\in\Theta_q}f_rf_s
 =\sum_{r\in\Theta_p\bigcap\Theta_q}f_r.
$$
Hence $\Theta_p$ and $\Theta_q$ are disjoint. On the other hand,
$$
 \sum_{p=1}^m f_p
 =\sum_{p=1}^m \theta(f_p)
 =\sum_{p=1}^m\sum_{q\in\Theta_p} f_q
 =\sum_{q\in \bigcup_{p=1}^m \Theta_p} f_q.
$$
Thus $\bigcup_{p=1}^m \Theta_p=\{1,\ldots,m\}$.

The partition $(\Theta_p)_{p=1,\ldots ,m}$ defines a unique function $u:\{1,\ldots,m\}\rightarrow\{1,\ldots,m\}$ given by $u(q)=p$ for all $q\in \Theta_p$. Hence, for an arbitrary $p\in\{1,\ldots,m\}$, we have
\begin{equation}\label{eq:theta}
\theta(f_p)=\sum_{q=1}^m \delta_{p,u(q)}f_q .
\end{equation}

\begin{lemma}\label{le:kernel}
The algebra endomorphism $\theta$ is idempotent if, and only if, the function $u$ is so. Moreover, the kernel of $\theta$ is given by
$$
 \Ker\theta=\langle f_p \mid p\notin \Im u \rangle.
$$
\end{lemma}
\begin{proof}
For $p\in \{1,\ldots,m\}$, we have $\theta^2(f_p)=\theta(f_p)$. By relation \eqref{eq:theta} we get
\begin{equation}\label{eq:u}
 \delta_{p,u^2(r)}=\delta_{p,u(r)},\mbox{ for all }p,r\in\{1,\ldots,m\}.
\end{equation}
Conversely, if \eqref{eq:u} holds, then $\theta$ is idempotent. We deduce that $\theta$ is an idempotent if, and only if, $u$ is so.

The partition associated to $\theta$ is given by $\Theta_p=u^{-1}(p)$. Thus $\Ker\theta$ is generated, as a vector space, by all $f_p$ such that $u^{-1}(p)=\emptyset$.
\end{proof}

\begin{theorem}\label{te:rank1}
Let $A=K^m$ and  let $\Gamma$ be a quiver of reduced rank one without cycles of length two. A splitted, unital and factorizable representation of $\Gamma$ is uniquely defined by a set of idempotent functions $u_i:\{1,\ldots,m\}\rightarrow\{1,\ldots,m\}$ with $i\in\Gamma^0$, satisfying the following condition: if
$\alpha\in\Gamma^1$ is not a loop and $p\in\{1,\ldots,m\}$, then $u_{s(\alpha)}(p)=p$ or $u_{t(\alpha)}(p)=p$.
\end{theorem}
\begin{proof}
By Theorem~\ref{te:Rank1Gen}, a splitted, unital and factorizable representation of $\Gamma$ is given by a family $(\varphi_i)_{i\in\Gamma^0}$ of idempotent algebra endomorphisms of $A=K^m$ such that, for any arrow $\alpha\in\Gamma^1$ that is not a loop, we have
$$
 \Ker\varphi_{s(\alpha)}\Ker\varphi_{t(\alpha)}=0.
$$
By Lemma~\ref{le:kernel}, for every $i\in\Gamma^0$ the algebra map $\varphi_i$ corresponds to an idempotent function $u_i:\{1,\ldots,m\}\rightarrow\{1,\ldots, m\}$. Let $\alpha$ be an arrow that is not a loop. Since
\[
 \Ker\varphi_i=\langle f_p\mid p\notin\Im\, u_i\rangle
\]
we deduce that
\begin{align*}
 \Ker\varphi_{t(\alpha)}\Ker\varphi_{s(\alpha)}
 &=\langle f_pf_q\mid p\notin\Im u_{s(\alpha)}, q\notin\Im u_{t(\alpha)}\rangle\\
 &=\langle f_p \mid p\notin \Im u_{s(\alpha)}\cup\Im u_{t(\alpha)}\rangle.
\end{align*}
Therefore $\Ker\varphi_{t(\alpha)}\Ker\varphi_{s(\alpha)}=0$ if, and only if, $\Im u_{s(\alpha)}\cup\Im u_{t(\alpha)}=\{1,\ldots,m\}$. Thus, for any $p\in\{1,\ldots,m\}$, $p\in\Im u_{s(\alpha)}$ or $p\in\Im u_{t(\alpha)}$. Since $u_i$ is idempotent for every $i\in\Gamma^0$, it follows that $\Im u_i=\{p\mid\;u_i(p)=p\}$, so the theorem is now proved.
\end{proof}

\section{Absolutely reducible representations of an algebra}\label{sec:absolutelyreducible}

We fix a finite-dimensional algebra $A$ over a field $K$ and a vector space $W$ of dimension $n$. Our aim in this section is to find all $A$-module structures on $W$ which are absolutely reducible in the following sense:

\begin{definition}
An $A$-module $W$ is called \dtext{absolutely reducible} (or \dtext{diagonalizable}) if there are $n$ submodules $W_1$, \ldots, $W_n$ of dimension one such that $W=\bigoplus_{i=1}^nW_i$.
\end{definition}

Let us fix a basis $\{w_1,\ldots,w_n\}$ on $W$. Thus a module structure on $W$ is given by some $K$-linear maps
$\omega_{ij}:A\longrightarrow{K}$, uniquely defined such that
\begin{equation}\label{eq:module}
 a{w_i}=\sum_{j=1}^n\omega_{ji}(a)w_j,\qquad\forall{i=1,\ldots,n}\mbox{\quad and\quad}\forall{a\in{A}}.
\end{equation}
The maps $\omega_{ij}$ are subject to the following relations:
\begin{align}
 &\sum_{k=1}^n\omega_{ik}(a)\omega_{kj}(b)=\omega_{ij}(ab).\label{eq:mod1}\\
 &\omega_{ij}(1_A)=\delta_{i,j}\label{eq:mod2}.
\end{align}
It is not difficult to see that, conversely, maps $\omega_{ij}$ satisfying \eqref{eq:mod1} and \eqref{eq:mod2} define a module structure on $W$.

\begin{example}\label{ex:module}
Let $\tau:K^n\otimes{A}\longrightarrow{A\otimes{K^n}}$ be a twisting  map, where $A=K^m$. Let
$\boldsymbol{E}_\tau=\{E_{ij}\}_{i,j=1,\ldots,n}$ be the corresponding set of $K$-linear endomorphisms of $A$. We fix a basis
$\{f_1,\ldots,f_m\}$ on $K^m$ and write
\[
 E_{ij}(a)=\sum_{p=1}^mE_{ij}^p(a)f_p,\qquad\forall{i,j=1,\ldots,n}\mbox{\quad and\quad}\forall{a\in{A}}.
\]
Hence, the set $\{E_{ij}^p\}_{i,j,=1,\ldots,n}$ satisfies \eqref{eq:mod1} and \eqref{eq:mod2}, for any $p=1,\ldots,m$.
\end{example}

We now assume that $W$ is absolutely reducible. Let $W=W_1\oplus\cdots\oplus{W_n}$ be the corresponding decomposition.
For each $i=1,\ldots,n$, we may choose a non-zero element $w_i'\in W_i$. Since $W_i$ is an one-dimensional submodule of $W$, there is an algebra map $\alpha_i:A\longrightarrow{K}$ such that
\[
 a{w_i'}=\alpha_i(a)w_i',\qquad\forall{a\in{A}}.
\]
As $\{w_1',\ldots,w_n'\}$ is another basis on $W$, there is an invertible matrix $X\in{GL_n(K)}$ such that
$X=(a_{ij})_{i,j=1,\ldots,n}$ and
\[
 w_i'=\sum_{j=1}^na_{ji}w_j.
\]
If $\{\omega_{ij}\}_{i,j=1,\ldots,n}$ are the maps that define the module structure on $W$, then we get
\[
 (\omega_{ij}(a))_{i,j=1,\ldots,n}=X\cdot
 \left(\begin{array}{ccc}
 \alpha_1(a)&\cdots&0\\\vdots&\ddots&\vdots\\0&\cdots&\alpha_n(a)
 \end{array}\right)
 \cdot{X^{-1}},\qquad\forall{a\in{A}}.
\]
For simplicity, we shall denote the matrix $(\omega_{ij})_{i,j=1,\ldots,n}$ by $\boldsymbol{\omega}$. Note
that the elements of $\boldsymbol{\omega}$ are in $\Hom_K(A,K)$. Hence the above equality can be written as
\begin{equation}\label{eq:phi}
 \boldsymbol{\omega}=X\cdot
 \left(\begin{array}{ccc}
 \alpha_1&\cdots&0\\\vdots&\ddots&\vdots\\0&\cdots&\alpha_n
 \end{array}\right)
 \cdot{X^{-1}}.
\end{equation}
Thus we have just proved the following Lemma.

\begin{lemma}\label{le:representations}
Let $A$ be a $K$-algebra. For every module structure on $W$ which is absolutely reducible, there are $X\in{GL_n(K)}$ and $\alpha_1$, \ldots, $\alpha_n\in\Alg_K(A,K)$ such that the matrix $\boldsymbol{\omega}$ defining the action of $A$ on $W$ satisfies \eqref{eq:phi}.
\end{lemma}

In view of this lemma, the set of characters of $A$ plays an important r\^{o}le in the description of absolutely reducible representations of an algebra $A$. By Dedekind's theorem on linear independence of characters,
%
%
any set of distinct algebra morphisms from $A$ to $K$ is linearly independent. In particular
%
%
$\Alg_K(A,K)$ is a finite set. We shall denote it by
\begin{equation}\label{eq:Alg}
 \Alg_K(A,K)=\{\theta_1,\ldots,\theta_r\}.
\end{equation}
%

\begin{remark}
Lemma~\ref{le:representations} can be rephrased as follows: if $W$ is an absolutely reducible module then there are
$X\in{GL_n(K)}$ and a set map $u:\{1,\ldots,n\}\longrightarrow\{1,\ldots,r\}$ such that
\begin{equation}\label{eq:phi1}
 \boldsymbol{\omega}=X\cdot\left(\begin{array}{ccc}
 \theta_{u(1)}&\cdots&0\\\vdots&\ddots&\vdots\\0&\cdots&\theta_{u(n)}
 \end{array}\right)\cdot{X^{-1}},
\end{equation}
where $\boldsymbol{\omega}$ is the set of $K$-linear maps associated to $W$. We shall say that $W$ is defined by the matrix $X$ and the map $u$. Of course, $X$ and $u$ are not uniquely determined by $W$. Given such a map $u$, our aim is to find a matrix $X_0\in{GL_n(K)}$ such that, together with $u$, it defines $W$ and has as many zero elements as possible.

Throughout the remaining of this section we fix an $n$-dimensional absolutely reducible representation $W$. Let $\boldsymbol{\omega}$ denote the corresponding set of $K$-linear maps and let $u:\{1,\ldots,n\}\longrightarrow\{1,\ldots,r\}$ be a function that defines $W$. Let
\[
 \Im(u)=\{i_1,\dots,i_s\},
\]
where $1\leq{i_1}<{i_2}<\ldots<{i_s}\leq{r}$. We denote $u^{-1}(i_k)$, the fiber of $u$ over $i_k$, by $F_k$. Hence $\mathcal{F}=\{F_1,\ldots,F_s\}$ is a partition of $\{1,\ldots,n\}$.
\end{remark}

\begin{definition}
Let $\mathcal{F}=\{F_1,\ldots,F_s\}$ be a partition of $\{1,\ldots,n\}$. For every matrix $X\in{M_n(K)}$, we define the
$\mathcal{F}$-\dtext{blocks} of $X$ by
\[
 X^{ij}=(x_{p,q})_{p\in{F_i},q\in{F_j}},
\]
where $x_{p,q}$ are the elements of $X$ and $i,j$ are in $\{1,\ldots,s\}$.
\end{definition}

\begin{example}\label{ex:blocks}
Let $u:\{1,2,3,4\}\longrightarrow\{1,2,3,4\}$ given by $u(1)=4$, $u(2)=u(4)=2$, $u(3)=1$. Hence $F_1=\{3\}$, $F_2=\{2,4\}$ and $F_3=\{1\}$. Thus, for any $X\in{M_4(K)}$ there are nine $\mathcal{F}$-blocks $(X^{kl})_{k,l=1,2,3}$. For example
\[
 X^{21}=\left(\begin{array}{cc}x_{23}\\x_{43}\end{array}\right),\quad
 X^{22}=\left(\begin{array}{cc}x_{22}&x_{24}\\x_{42}&x_{44}\end{array}\right),\quad
 X^{32}=\left(\begin{array}{cc}x_{12}&x_{14}\end{array}\right).\]
\end{example}

\begin{lemma}\label{le:H}
Let $X\in{GL_n(K)}$, $u:\{1,\ldots,n\}\to\{1,\ldots,r\}$ be a set map and $\mathcal{F}$ be the set of fibers of $u$.
\begin{enumerate}[(1)]
\item
The set
\[
 H_u=\{Y\in{GL_n(K)}\mid\;Y^{kl}=0,\mbox{ for }k\neq{l}\}
\]
is a subgroup of $GL_n(K)$.
\item
If $X\in{GL_n(K)}$ and $Y\in{H_u}$, then $X$ and $XY$ define, together with $u$, the
same representation. Moreover,
\[
 (XY)^{kl}=X^{kl}Y^{ll},\qquad\forall{k,l\in\{1,\ldots,s\}}.
\]
\item
Let $X$, $Y\in{GL_n(K)}$. If $X$ and $Y$ define the same representation, then there is
$Z\in{H_u}$ such that $Y=XZ$.
\end{enumerate}
\end{lemma}
\begin{proof}
(1) %
Let us show that $H_u$ contains all matrices in $GL_n(K)$ that commute with
\[
 \boldsymbol{\theta}=
 \left(\begin{array}{ccc}\theta_{u(1)}&\cdots&0\\\vdots&\ddots&\vdots\\0&\cdots&\theta_{u(n)}\end{array}\right).
\]
Indeed, for $Y\in{GL_n(K)}$, we have $Y\boldsymbol{\theta}=\boldsymbol{\theta}Y$ if, and only if,
\[
 y_{ij}(\theta_{u(i)}-\theta_{u(j)})=0,\qquad\forall{i,j=1,\ldots,n}.
\]
These equalities are equivalent to $y_{ij}=0$ for any pair $(i,j)$ such that $i$ and $j$ are not in the same fiber of $u$. In conclusion, $Y$ commutes with $\boldsymbol{\theta}$ if, and only if, $Y\in{H_u}$. Since the set of commuting matrices with $\boldsymbol{\theta}$ is a subgroup of $GL_n(K)$, the first part of the lemma is proved.

(2) %
Let $Y\in{H_u}$. Since $Y\boldsymbol{\theta}Y^{-1}=\boldsymbol{\theta}$, then obviously $XY$ and $Y$ define the same representation. It remains to prove the formula that describes the $\mathcal{F}$-blocks of $XY$. Let $k$, and $l$ be arbitrary in $\{1,\ldots,s\}$. We choose $i\in{F_k}$ and $j\in{F_l}$. If $z_{ij}$ is the element of $XY$ in
the spot $(i,j)$, then
\[
 z_{ij}=\sum_{p=1}^nx_{ip}x_{pj}=\sum_{l'=1}^s\sum_{p\in{F_{l'}}}x_{ip}y_{pj}.
\]
By definition of $H_u$, for $l'\neq{l}$ we have $y_{pj}=0$, as $y_{pj}$ is an element in $Y^{l'l}$. So
\[
 z_{ij}=\sum_{p\in{F_l}}x_{ip}y_{pj}.
\]
Obviously the element in the  $(i,j)$-spot of $X^{kl}Y^{ll}$ is $\sum_{p\in{F_l}}x_{ip}x_{pj}$, so this part of the lemma is also proved.

(3) If $X$ and $Y$, together with $u$, define the same representation, then $Z=X^{-1}Y$ commutes with $\boldsymbol{\theta}$. Then $Z\in{H_u}$ and $Y=XZ$.
\end{proof}

\begin{theorem}
We keep the notation from Lemma~\ref{le:H}. Then
\[
 H_u\cong{GL_{n_1}(K)}\times\cdots\times{GL_{n_s}(K)},
\]
where $n_k=\#{F_k}$ for $k=1,\ldots,s$. The right action of $H_u$ on $GL_n(K)$ given by right matrix multiplication is defined on $\mathcal{F}$-blocks by the formula
\begin{equation}\label{eq:30}
 (XY)^{kl}=X^{kl}Y^{ll},\qquad\forall{k,l=1,\ldots,s}.
\end{equation}
for any $X\in GL_n(K)$ and $Y\in H_u$.
\end{theorem}
\begin{proof}
Let $Y\in{H_u}$. By Lemma~\ref{le:H}(1) the $\mathcal{F}$-blocks $Y^{kl}$ are trivial, for any $k\neq{l}$. Hence the map $\Phi:H_u\longrightarrow{GL_{n_1}(K)}\times\cdots\times{GL_{n_s}(K)}$ given by  \[\Phi(Y):=(Y^{11},\ldots,Y^{ss})\] is well-defined. Note that $Y^{kk}$ is invertible, for any $k=1,\ldots,s$. Indeed, up to a permutation of the rows and columns of $Y$, we have
\[
Y=\left(\begin{array}{ccc}Y_{11}&\cdots&0\\\vdots&\ddots&\vdots\\0&\cdots&Y_{ss}\end{array}\right).
\]
By the formula in Lemma~\ref{le:H}(2), it follows that $\Phi$ is an isomorphism of groups. By the same formula, we can describe the action of $H_u$ on $GL_n(K)$ as in \eqref{eq:30}.
\end{proof}

Recall that our aim is to find a ``{normalized}" matrix $X\in{GL_n(K)}$ that defines, together with $u$, a given
absolutely reducible representation $W$. We show that the action of $H_u$ on $M_n(K)$ can be splitted  in $s$ actions $M_{n,n_k}(K)\times GL_{n_k}(K)\to M_{n,n_k}(K)$, for $k=1,\ldots,s$.

Let $X=(x_{ij})_{i,j=1,\dots,n}$ and let $u:\{1,\ldots,n\}\longrightarrow\{1,\ldots,r\}$. Recall that
$\mathcal{F}$ denotes the partition $\{F_1,\dots,F_s\}$ associated to $u$. For $1\leq k\leq s$, we define $X^k\in M_{n,n_k}(K)$ to be the matrix whose elements are $x_{ij}$, where $i\in\{1,\ldots,n\}$ and $j\in{F_k}$.

Each $X^k$ is made of the $\mathcal{F}$-blocks $X^{1k}$, $X^{2k}$, \ldots, $X^{sk}$. Hence we have
\begin{equation}\label{eq:action}
 (XY)^k=X^kY^{kk}
\end{equation}
for any $Y\in H_u$.

 The group $GL_{n_k}(K)$ acts on $M_{n,n_k}(K)$
by right multiplication. Thus, the action of $H_u$ on $GL_{n}(K)$
can be recovered from these $s$ actions by the relation
\eqref{eq:action}.

\begin{definition}
The symmetric group $S_n$ acts on $M_{n,m}(K)$ by row permutation,
$(\sigma,X)\mapsto{X_\sigma}$, where
\[
 X_\sigma:=\left(\begin{array}{ccc}
 x_{\sigma(1)1}&\cdots&x_{\sigma(1)m}\\
 x_{\sigma(2)1}&\cdots&x_{\sigma(2)m}\\
 \vdots&&\vdots\\
 x_{\sigma(n)1}&\cdots&x_{\sigma(n)m}\\
 \end{array}\right)
\]
\end{definition}

\begin{theorem}\label{th:normalizedX}
Let $X\in{GL_n(K)}$ and let
$u:\{1,\ldots,n\}\longrightarrow\{1,\ldots,r\}$ be an arbitrary set
map. The orbit of $X$ with respect to the action of $H_u$ on
$GL_n(K)$ contains a matrix $\overline{X}$ such that
\[
 \overline{X}^k=\left(\begin{array}{c}I_{n_k}\\Z_k\end{array}\right)_{\sigma_k},
\]
where $Z_k$ is a certain matrix in $M_{n-n_k,n_k}(K)$ and $\sigma_k\in{S_n}$ for any $k\in\{1,\ldots,s\}$.
\end{theorem}
\begin{proof}
Since $X\in{GL_n(K)}$, the blocks $X^1$, \ldots, $X^s$ are matrices
of rank $n_1$, \ldots, $n_s$ respectively (the columns of $X$ are
linearly independent, so the columns of each $X^k$ are so). Thus
there are $n_k$ rows of $X^k$ that are linearly independent. There
is $\tau_k\in{S_n}$ such that $X_{\tau_k}^k$ has the first $n_k$
rows linearly independent. It follows that
\[
 X_{\tau_k}^k=\left(\begin{array}{c}Y_{k}\\X'_k\end{array}\right)
\]
with $Y_k\in{GL_{n_k}(K)}$ and $X'_k\in{M_{n-n_k,n_k}(K)}$. Hence
\[
 X^k=\left(\begin{array}{c}I_{n_k}\\X'_kY^{-1}_k\end{array}\right)_{\tau^{-1}_k}\cdot{Y_k},
\]
so we may take $Z_k=X'_kY^{-1}_k$ and $\sigma_k=\tau^{-1}_k$.
\end{proof}

\begin{remark}
The matrices $Z_1$, \ldots, $Z_s$ and the permutations $\sigma_1$, \ldots, $\sigma_s$ are not uniquely determined.
\end{remark}

\begin{definition}
The matrix $\overline{X}$ as in Theorem~\ref{th:normalizedX} will be called a \emph{normalized invertible matrix}.
\end{definition}

\begin{corollary}
Let $W$ be an absolutely reducible representation of $A$ of
dimension $n$. Then there are a normalized invertible matrix
$X\in{GL_n(K)}$ and a set map
$u:\{1,\ldots,n\}\longrightarrow\{1,\ldots,r\}$ that define the
action of $A$ on $W$.
\end{corollary}
\begin{proof}
Any absolutely reducible representation is defined by a matrix $X_0$
and a map $u:\{1,\ldots,n\}\longrightarrow\{1,\ldots,r\}$. By the
Theorem~\ref{th:normalizedX}, there is a normalized invertible
matrix $X$ in the orbit of $X_0$, with respect to the action of
$H_u$ on $GL_n(K)$. We have already noticed that two matrices in the
same orbit define the same representation. So $W$ is defined by $X$
and $u$.
\end{proof}

\begin{example}
We keep the notation from Example~\ref{ex:blocks}. An arbitrary element in $H_u$ has the following form
\[
 Y=\left(\begin{array}{cccc}
 y_{11}&0&0&0\\
 0&y_{22}&0&y_{24}\\
 0&0&y_{33}&0\\
 0&y_{42}&0&y_{44}
 \end{array}\right),
\]
where $y_{11}$ and $y_{33}$ are non-zero and
$y_{22}\,y_{44}\neq{y_{24}}\,y_{42}$. For an  normalized invertible
matrix $X\in{GL_4(K)}$ we have
\[
 X^1=\left(\begin{array}{cc}1\\x_{21}\\x_{31}\\x_{41}\end{array}\right)_{\sigma_1},\quad
 X^2=\left(\begin{array}{cc}1&0\\0&1\\x_{32}&x_{34}\\x_{42}&x_{44}\end{array}\right)_{\sigma_2},\quad
 X^3=\left(\begin{array}{cc}1\\x_{23}\\x_{33}\\x_{43}\end{array}\right)_{\sigma_3}.
\]
The permutations  $\sigma_1$, $\sigma_2$, $\sigma_3\in{S_4}$ have to be chosen such that $\det(X)\neq0$. For example
\[
 X=\left(\begin{array}{cccc}
 x_{21}&0&x_{43}&1\\
 x_{41}&x_{42}&1&x_{44}\\
 1&x_{32}&x_{23}&x_{34}\\
 x_{21}&1&x_{33}&0
 \end{array}\right)
\]
is a normalized invertible matrix if, and only if, $\det(X)\neq0$.
\end{example}

We end this section by analyzing in detail the case of
two-dimensional absolutely reducible representations.

\begin{proposition}\label{pr:normalized}
Let $X\in{GL_2(K)}$ and let $u:\{1,2\}\longrightarrow\{1,\ldots,r\}$
be a set map. Then there are $x$, $y\in{K}$, with $xy\neq1$ such
that the orbit of $X$ with respect to the action of $H_u$ on
$GL_2(K)$ contains one of the matrices
\[
 X_1=\left(\begin{array}{cc}1&x\\y&1\end{array}\right)\quad\mbox{ or }\quad
 X_2=\left(\begin{array}{cc}x&1\\1&y\end{array}\right).
\]
\end{proposition}
\begin{proof}
We start by considering the case $u(1)=u(2)=i_1\in\{1,\ldots,r\}$.
Hence the partition $\mathcal{F}$ associated to $u$ has one set
$F_1=\{1,2\}$ and $H_u=GL_2(K)$. In conclusion, $X$ and
$I_2=XX^{-1}$ are in the same orbit, so in this case we can take
$x=y=0$ in $X_1$.

Let now take $u:\{1,2\}\longrightarrow\{1,\ldots,r\}$ such that
$u(1)\neq{u(2)}$. If $\Im(u)=\{i_1,i_2\}$, with
$1\leq{i_1}<i_2\leq{r}$, then the partition $\mathcal{F}$ is given
either by $F_1=\{1\}$ and $F_2=\{2\}$ or $F_1=\{2\}$ and
$F_2=\{1\}$. In both cases we get that a normalized invertible
matrix is of the following type:
\[
 Z_1=\left(\begin{array}{cc}1&1\\x_1&y_1\end{array}\right),\;
 Z_2=\left(\begin{array}{cc}x_2&y_2\\1&1\end{array}\right),\;
 Z_3=\left(\begin{array}{cc}1&x_3\\y_3&1\end{array}\right),\;
 Z_4=\left(\begin{array}{cc}x_4&1\\1&y_4\end{array}\right).
\]
Since $Z_i$ is invertible we have $x_i\neq{y_i}$ for
$i\in\{1,2\}$, and $x_iy_i\neq1$ for $i\in\{3,4\}$. We know that
in the orbit of $X$ there is a matrix $Z_i$, for a certain
$i\in\{1,2,3,4\}$. If either $i=3$ or $i=4$, one can take
$X_1=Z_3$ or $X_2=Z_4$.

Let us assume that $i=1$. Note that
\[
 H_u=\left\{\left(\begin{array}{cc}y_1&0\\0&y_2\end{array}\right)\mid\;y_1y_2\neq0\right\}
\]
and either $x_1\neq0$ or $y_1\neq0$. If $x_1\neq0$ we have
\[
 Z_1\cdot
 \left(\begin{array}{cc}x_1^{-1}&0\\0&1\end{array}\right)=
 \left(\begin{array}{cc}x_1^{-1}&1\\1&y_1\end{array}\right),
\]
so in the orbit of $Z_1$ (and hence of $X$) there is a matrix of
type $\left(\begin{array}{cc}x&1\\1&y\end{array}\right)$. In the
case when $y_1\neq0$ one can show that the orbit of $X$ contains a
matrix of the same type. In the case $i=2$ we can proceed
similarly to show that the orbit of $X$ contains a matrix $X_2$ as
in the statement of the theorem.
\end{proof}

Let $W$ be a vector space of dimension 2. We fix a basis
$\{w_1,w_2\}$ on úWú. We want to classify all $A$-module
structures on $W$ which are absolutely reducible.

\begin{theorem}\label{th_rep2}
Let $A$ be a $K$-algebra and let $W$ be a two-dimensional
absolutely reducible $A$-module. Then there are $x$, $y\in{K}$ and
$\alpha_1$, $\alpha_2\in\Alg_K(A,K)$ such that $xy\neq1$ and for
any $a\in{A}$
\begin{align}
 a{w_1}&=\frac{1}{1-xy}\left[\alpha_1(a)-xy\alpha_2(a)\right]w_1
 +\frac{y}{1-xy}\left[\alpha_1(a)-\alpha_2(a)\right]w_2,\label{eq:W1}\\
 a{w_2}&=\frac{-x}{1-xy}\left[\alpha_1(a)-\alpha_2(a)\right]w_1
 +\frac{1}{1-xy}\left[-xy\alpha_1(a)+\alpha_2(a)\right]w_2,\label{eq:W2}
\end{align}
\end{theorem}
\begin{proof}
The representation $W$ is defined by a function
$u:\{1,2\}\longrightarrow\{1,\ldots,r\}$ and a normalized invertible
matrix $X$. Since two matrices in the same orbit with respect to the
action of $H_u$ define the same representation we may assume, in
view of the previous proposition, that
\[
 X=\left(\begin{array}{cc}1&x\\y&1\end{array}\right)\quad\mbox{ or }\quad
 X=\left(\begin{array}{cc}x&1\\1&y\end{array}\right).
\]
Let $\boldsymbol{\omega}:=(\omega_{ij})_{i,j=1,2}$ be the matrix
associated to the $A$-module $W$ as in \eqref{eq:module}. There
are $\theta_{u(1)}$ and $\theta_{u(2)}$ in $\Alg_K(A,K)$ such that
\[
 \boldsymbol{\omega}=X\left(\begin{array}{cc}\theta_{u(1)}&0\\0&\theta_{u(2)}\end{array}\right)X^{-1}.
\]
We first consider the case
$X=\left(\begin{array}{cc}1&x\\y&1\end{array}\right)$. If
$\alpha_1:=\theta_{u(1)}$ and $\alpha_2:=\theta_{u(2)}$, then a
straightforward computation shows us that
\begin{equation}\label{eq:omega}
 \boldsymbol{\omega}=\left(\begin{array}{cc}
 {(1-xy)^{-1}}\left[\alpha_1-xy\alpha_2\right]&-x(1-xy)^{-1}\left[\alpha_1-\alpha_2\right]\\
y(1-xy)^{-1}\left[\alpha_1-\alpha_2\right]&(1-xy)^{-1}\left[-xy\alpha_1+\alpha_2\right]\\
 \end{array}\right)
\end{equation}
It is easy to see that the module structure defined by $\boldsymbol{\omega}$ in this case is as in \eqref{eq:W1} and \eqref{eq:W2}.

If $X=\left(\begin{array}{cc}x&1\\1&y\end{array}\right)$, then, for
$\alpha_1=\theta_{u(2)}$ and $\alpha_2=\theta_{u(1)}$, one can show
that the corresponding action  also satisfies relations
\eqref{eq:W1} and \eqref{eq:W2}.
\end{proof}

\section{Splitted, unital and factorizable representations of cycles of length 2}\label{sec:cyclesoflength2}

In this section, for $A=K^m$, we shall classify all splitted, unital
and factorizable representations of a cycle $\Gamma$ of length 2. We
use the standard quiver-notation $\Gamma^0=\{1,2\}$ and
$\Gamma^1=\{\alpha_1,\alpha_2\}$ in order to represent the quiver
below.
$$
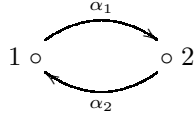
\xy
 (0,0)*+{1\;\circ}="a",
 (20,0)*+{\circ\;2}="b",
 \ar@/^3ex/ "a";"b"^{\alpha_1}
 \ar@/_-3ex/ "b";"a"^{\alpha_2}
\endxy
$$
A representation of $\Gamma$ is defined by a set
$\{\varphi_1,\varphi_2,\varphi_{\alpha_1},\varphi_{\alpha_2}\}$ of
$K$-linear maps. For $i=1,2$, there are
$\varphi_i^1,\ldots,\varphi_i^m,\varphi_{\alpha_{i}}^1,\ldots,\varphi_{\alpha_{i}}^m\in\Hom_K(A,K)$
such that
\begin{align}
 \varphi_i(a)&=\sum_{p=1}^m\varphi_i^p(a)f_p,\qquad\forall{a\in{A}},
\label{fi_i}\\
 \varphi_{\alpha_i}(a)&=\sum_{p=1}^m\varphi_{\alpha_i}^p(a)f_p,\qquad\forall{a\in{A}}.\label{fi_ai}
\end{align}
Here $\{f_1,\ldots,f_m\}$ denotes the canonical basis on $K^m$. If $\{f_1^*,\ldots,f_m^*\}$ is the dual basis of
$\{f_1,\ldots,f_m\}$, then
\[
 \Alg_K(K^m,K)=\{f_1^*,\ldots,f_m^*\}.
\]
Let us assume that
$\{\varphi_1,\varphi_2,\varphi_{\alpha_1},\varphi_{\alpha_2}\}$
defines a splitted, unital and factorizable representation of
$\Gamma$. Since it is unital and factorizable, it is not difficult
to see that the matrix
$\boldsymbol{\omega^p}=(\omega_{ij}^p)_{i,j=1,2}$ given by
\begin{equation}\label{eq:omega_p}
 \boldsymbol{\omega^p}=\left(\begin{array}{cc}\varphi_1^p&\varphi_{\alpha_1}^p\\\varphi_{\alpha_2}^p&\varphi_2^p\end{array}\right)\qquad
 \forall{p\in\{1,\ldots,m\}}
\end{equation}
defines an $A$-module structure on a vector space $W$ of dimension 2. If we fix a basis $\{w_1,w_2\}$ on $W$, the corresponding module structure is defined by
\[
 a{w_i}=\sum_{j=1}^n\omega_{ji}^p(a)w_j,\qquad\forall{a\in{A}},\;\forall{i\in\{1,2\}}.
\]
To emphasize the fact that this module structure depends on $p\in\{1,\ldots,m\}$, we shall denote it by $W_p$.

\begin{theorem}\label{th:duplicates}
Let $A=K^m$ and let $\Gamma$ be a cycle of length 2. If $\boldsymbol{\omega^p}$ is defined as in \eqref{eq:omega_p}, then there are $a_1$, \ldots, $a_m\in{K}$ and $u:\{1,\ldots,m\}\longrightarrow\{1,\ldots,m\}$ such that
\begin{equation}\label{eq:op}
 \boldsymbol{\omega^p}=\left(\begin{array}{cc}
 a_pf_p^*+(1-a_p)f_{u(p)}^*&a_p(f_p^*-f_{u(p)}^*)\\
 (1-a_p)(f_p^*-f_{u(p)}^*)&(1-a_p)f_p^*+a_pf_{u(p)}^*
 \end{array}\right).
\end{equation}
\end{theorem}
\begin{proof}
Since $A=K^m$ is a semisimple algebra and every simple $A$-module is
one dimensional, it follows that every representation of $A$ is
absolutely reducible. In particular, $W_p$ is so, for any
$p\in\{1,\ldots,m\}$. By the proof of Theorem~\ref{th_rep2}, there
are $x^p$, $y^p\in{K}$ with $x^py^p\neq1$, and $\alpha_1^p$,
$\alpha_2^p\in\Alg_K(K^m,K)$ such that
\[
 \boldsymbol{\omega^p}=\left(\begin{array}{cc}
 {(1-x^py^p)^{-1}}\left[\alpha_1^p-x^py^p\alpha_2^p\right]
 &{-x^p}{(1-x^py^p)^{-1}}\left[\alpha_1^p-\alpha_2^p\right]\\
 {y^p}{(1-x^py^p)^{-1}}\left[\alpha_1^p-\alpha_2^p\right]
 &{(1-x^py^p)^{-1}}\left[-x^py^p\alpha_1^p+\alpha_2^p\right]
 \end{array}\right).
\]
On the other hand, since $\{\varphi_1,\varphi_2,\varphi_{\alpha_1},\varphi_{\alpha_2}\}$ defines a splitted representation of $\Gamma$, we have
$\varphi_1+\varphi_{\alpha_2}=\Id_A$ and $\varphi_2+\varphi_{\alpha_1}=\Id_A$. These
relations are equivalent to
\begin{equation}\label{eq:split2}
 \varphi_1^p+\varphi_{\alpha_2}^p=f_p^*\quad\mbox{ and }\quad
 \varphi_2^p+\varphi_{\alpha_1}^p=f_p^* \qquad \text{for any
 $p\in\{1,\ldots, n\}$},
\end{equation}
therefore the sum of the elements of each column of
$\boldsymbol{\omega^p}$ is equal to $f_p^*$.

If $\alpha_1^p=\alpha_2^p$, then the above condition imposed on the elements of the column in $\boldsymbol{\omega^p}$ means that $\alpha_1^p=\alpha_2^p=f_p^*$. In this case we take $a_p=1$ and we define $u(p)=p$.

Now let us assume that $\alpha_1^p\neq\alpha_2^p$. Then these
algebra morphisms are linearly independent. By imposing that
$\boldsymbol{\omega^p}$ satisfies \eqref{eq:split2}, we get
\[
 \frac{y^p+1}{1-x^py^p}\alpha_1^p-\frac{y^p(x^p+1)}{1-x^py^p}\alpha_2^p
 =f_p^*
 =\frac{-x^p(y^p+1)}{1-x^py^p}\alpha_1^p+\frac{x^p+1}{1-x^py^p}\alpha_2^p.
\]
These equalities are possible only in the following two cases:

\emph{First case:} $y^p=-1$ and $\alpha_2^p=f_p^\ast$. As
$x^py^p\neq 1$, we have $x^p\neq-1$. Since $\alpha_1^p$ is an
algebra map, there is $u(p)\in\{1,\ldots,m\}$ such that
$\alpha_1^p=f_{u(p)}^*$. We now define $a^p:={x^p}{(1+x^p)^{-1}}$.
It is not difficult to see that $\boldsymbol{\omega^p}$ satisfies
\eqref{eq:op}.

\emph{Second case:} $x^p=-1$ and $\alpha_1^p=f_p^*$. As above, $y^p\neq-1$ and there is $1\leq u(p)\leq m$ such that
$\alpha_2^p=f_{u(p)}^*$. We conclude by defining $a^p:=(1+y^p)^{-1}$.
\end{proof}

Now we can classify all splitted, unital and factorizable
representations of a cycle of length 2.

\begin{theorem}\label{th:cycle2}
Let $A=K^m$  and let $\Gamma$ be a cycle of length 2. The set of
$K$-linear endomorphisms
$\boldsymbol{\varphi}:=\{\varphi_1,\varphi_2,\varphi_{\alpha_1},\varphi_{\alpha_2}\}$
defines a splitted, unital and factorizable representation of
$\Gamma$ if, and only if, there are $a_1$, \ldots, $a_m\in{K}$ and a
set map $u:\{1,\ldots,m\}\longrightarrow\{1,\ldots,m\}$ such that
\begin{align}
 \varphi_1&=\sum_{p=1}^m \left[a_pf_p^*+(1-a_p)f_{u(p)}^*\right]f_p,
 \quad
 \varphi_{\alpha_1}=\sum_{p=1}^m \left[a_p(f_p^*-f_{u(p)}^*)\right]f_p,\label{r1}\\
 \varphi_2&=\sum_{p=1}^m \left[(1-a_p)f_p^*+a_pf_{u(p)}^*\right]f_p,
 \quad
 \varphi_{\alpha_2}=\sum_{p=1}^m\left[(1-a_p)(f_p^*-f_{u(p)}^*)\right]f_p,\label{r2}
\end{align}
and, in addition, they satisfy the following conditions:
\begin{enumerate}[(1)]
\item
If $u(p)=p$, then $ a_p=0$.
\item
If $u(p)\neq p$, then $a_p+a_{u(p)}=1$. In addition, if $u^2(p)\neq p$, then $ a_p\in \{0,1\}$.
\end{enumerate}
\end{theorem}
\begin{proof}
In view of Theorem~\ref{th:duplicates}, there are $a_1$, \ldots,
$a_m\in{K}$ and $u$ such that the maps in $\boldsymbol{\varphi}$
satisfy the four equalities in the statement. Note that if $u(p)=p$,
then the $p$-th term of each sum does not depend on $a_p$.
Therefore, in this case, we may normalize $a_p$ by setting $a_p=0$.
Thus we may choose $a_1$, \ldots, $a_m$ such that the first
condition is verified.

Note that the above description of the maps $\varphi_1$,
$\varphi_2$, $\varphi_{\alpha_1}$ and $\varphi_{\alpha_2}$ was
obtained as a consequence of the fact that $\boldsymbol{\varphi}$
defines a unital and factorizable representation of $\Gamma$ such
that
$\varphi_1+\varphi_{\alpha_2}=\varphi_2+\varphi_{\alpha_1}=\Id_A$.
In conclusion, the representation associated to
$\boldsymbol{\varphi}$ is also splitted if, and only if,
$\varphi_{\alpha_1}$ and $\varphi_{\alpha_2}$ are idempotents.
Remark that $\{\varphi_1,\varphi_{\alpha_2}\}$ and
$\{\varphi_2,\varphi_{\alpha_1}\}$ then are complete sets of
orthogonal idempotents, as it is required.

Trivially, $\varphi_{\alpha_1}$ is an idempotent $K$-linear map if, and only if, $\varphi_{\alpha_1}^2(f_q)=\varphi_{\alpha_1}(f_q)$, for any
$q\in\{1,\ldots,m\}$. By definition of $\varphi_{\alpha_1}$, we get that these equalities hold if, and only if, for arbitrary
$p,q\in\{1,\ldots,m\}$,
\begin{equation}\label{eq:ap}
 a_p\left[(\delta_{p,q}-\delta_{u(p).q})a_p-(\delta_{u(p),q}-\delta_{u^2(p),q})a_{u(p)}\right]
 =a_p(\delta_{p,q}-\delta_{u(p),q}).
\end{equation}
Let us assume that $a_p\neq0$. By the normalization condition, then $u(p)\neq{p}$.

\emph{First case:} $u^2(p)=p$. Since $a_p\neq 0$, relation \eqref{eq:ap} is equivalent in this case to
\[
 (a_p+a_{u(p)}-1)(\delta_{p,q}-\delta_{u(p),q})=0,\qquad\forall{q\in\{1,\ldots,m\}}.
\]
Since $u(p)\neq{p}$, we deduce that $a_p+a_{u(p)}=1$.

\emph{Second case:} $u^2(p)\neq{p}$. Hence, for $q=p$, relation
\eqref{eq:ap} becomes $a_p=1$.  We want to prove that $a_{u(p)}=0$.
Note that if $u^2(p)=u(p)$, then this equality follows from the fact
that we have normalized the elements $a_1$, \ldots, $a_m$. So we may
assume that $u^2(p)$ and $u(p)$ are not equal. By taking $q=u(p)$ in
\eqref{eq:ap}, it follows that $a_{u(p)}=0$. Obviously,
$a_p+a_{u(p)}=1$ in this case too.

Summarizing, since $\varphi_{\alpha_1}$ is an idempotent $K$-linear
map then, for any $p$ such that $a_p\neq 0$, one of the following
conditions holds:
\begin{align}
 u^2(p)&=p, \text{ and then } a_p+a_{u(p)}=1;\label{eq01}\\
 u^2(p)&\neq p, \text{ and then } a_p=1 \text{ and } a_{u(p)}=0.\label{eq02}
\end{align}
Since $\varphi_{\alpha_2}=\sum_{r=1}^m
\left[(1-a_r)(f_r^*-f_{u(r)}^*)\right ]f_p$ is an idempotent too,
working with $1-a_p$ instead of $a_p$,  we can show that, for any
$p\in\{1,\ldots,m\}$ such that $a_p\neq 1$, one of the following two
conditions holds:
\begin{align}
 u^2(p)&=p, \text{ and then } a_p+a_{u(p)}=1;\label{eq11}\\
 u^2(p)&\neq p,\text{ and then } a_p=0 \text{ and } a_{u(p)}=1.\label{eq12}
\end{align}

Since either $a_p\neq 0$ or $a_p\neq 1$, the conditions \eqref{eq01}-\eqref{eq12} together imply that $a_1,\ldots,a_m$ and $u$ verify the second property in the statement of the theorem.

Now let us prove that, for any $a_1,\ldots,a_m$ and $u$ satisfying
(1) and (2), the set
$\boldsymbol{\varphi}:=\{\varphi_1,\varphi_2,\varphi_{\alpha_1},\varphi_{\alpha_2}\}$
defines a splitted, unital and factorizable representation of
$\Gamma$. By the proof of Theorem~\ref{th:duplicates},
$\boldsymbol{\varphi}$ defines a representation of $\Gamma$ which is
unital and factorizable. Of course, by construction,
$\varphi_1+\varphi_{\alpha_2}=\Id_A$ and
$\varphi_2+\varphi_{\alpha_1}=\Id_A$. By the proof of the other
implication, it results that conditions (1) and (2) imply that
$\varphi_1$ and $\varphi_2$ are idempotent maps. In conclusion, the
representation defined by $\boldsymbol{\varphi}$ is splitted too.
\end{proof}

\begin{corollary}
Let $\tau:K^n\otimes{K^m}\longrightarrow{K^m}\otimes{K^n}$ be a twisting map such that the connected components of $\Gamma$, the corresponding quiver, are all trivial excepting one which is a cycle of length 2. We label the vertices of $\Gamma$ as in the picture below.
$$
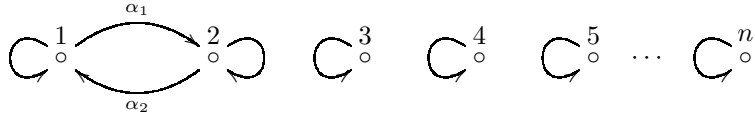
\xy
 (5,0)*+{\circ}="a",(5,3)*+{1},
 (25,0)*+{\circ}="b",(25,3)*+{2},
 (45,0)*+{\circ}="c", (45,3)*+{3},
 (60,0)*+{\circ}="d", (60,3)*+{4},
 (75,0)*+{\circ}="e", (75,3)*+{5},
 (82,0)*+{\cdots}, (95,0)*+{\circ}="f", (95,3)*+{n},
 \ar@/^3ex/ "a";"b"^{\alpha_1}
 \ar@/_-3ex/ "b";"a"^{\alpha_2}
\ar@(lu,ld) "c";"c"\ar@(lu,ld) "d";"d"\ar@(lu,ld)
"e";"e"\ar@(lu,ld) "f";"f"\ar@(lu,ld) "a";"a"\ar@(ru,rd) "b";"b"
\endxy
$$
There are $a_1$, \ldots, $a_m\in{K}$ and $u:\{1,\ldots,m\}\longrightarrow\{1,\ldots,m\}$ that satisfy
conditions (1) and (2) in Theorem~\ref{th:cycle2} such that
\begin{align*}
 \tau(e_1\otimes{x})&=\varphi_1(x)\otimes{e_1}+\varphi_{\alpha_1}(x)\otimes{e_2},\\
 \tau(e_2\otimes{x})&=\varphi_{\alpha_2}(x)\otimes{e_1}+\varphi_2(x)\otimes{e_2},\\
 \tau(e_i\otimes{x})&=x\otimes{e_i},\qquad\forall{i\geq3}.
\end{align*}
where $\varphi_1$, $\varphi_2$, $\varphi_{\alpha_1}$, $\varphi_{\alpha_2}$ are given by the formulae \eqref{r1} and \eqref{r2}.
\end{corollary}

\begin{remark}
If $n=2$, the corollary above gives the classification of
noncommutative duplicates of $K^m$ (cf. \cite{Cibils06a, Lopez08a}).
\end{remark}

We are going to classify all splitted, unital, and factorizable
representations of a connected quiver of rank 1 that contains a
cycle of length 2. We know that, removing the arrows of the cycle
from such a quiver, we get two rooted trees. These trees have their
roots in the vertices of the cycle and they are oriented. Note that
they might be trivial, that is they could have a unique vertex,
namely their root.

Let us denote the vertices of this cycle by $1$ and $2$, being
$3,4,\dotsc, n$ the other vertices of $\Gamma$, and denote by
$\alpha_1$ and $\alpha_2$ the arrows of the cycle, as in the picture
below.
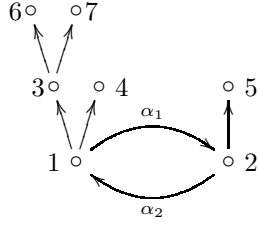
\[
\xy
 (0,0)*+{\circ}="a1",(20,0)*+{\circ}="a2",
 (-3,0)*+{1},(23,0)*+{2},
 (-3,10)*+{\circ}="b",(-5,10)*+{3},
 (3,10)*+{\circ}="c",(6,10)*+{4},(23,10)*+{5},
 (-6,20)*+{\circ}="d",(-8,20)*+{6},
 (0,20)*+{\circ}="f",(2,20)*+{7},
 (20,10)*+{\circ}="h",
 \ar@/^3ex/ "a1";"a2"^{\alpha_1}
 \ar@/^3ex/ "a2";"a1"^{\alpha_2}
 \ar "a1";"b"\ar "a1";"c"
 \ar "b";"d"
 \ar "b";"f"
 \ar "a2";"h"
\endxy
\]
We take a splitted, unital and factorizable representation of
$\Gamma$, which is defined by the $K$-linear maps
$\varphi_1,\dotsc,\varphi_n$ and $\{\varphi_{\alpha}\}$, where
$\alpha$ runs over the set of arrows which are not a loop. Note that
for $i\geq 3$, $\varphi_i$ is an idempotent algebra map, and thus,
there are idempotent maps $u_i:\{1,\dotsc,m\}\to \{1,\dotsc,m\}$
such that
\[
 \varphi_i(f_p) = \sum_{q=1}^m \delta_{p,u_i(q)}f_q,\ \forall\, p=1,\dotsc,m.
\]
Obviously,
$\boldsymbol{\varphi}=\{\varphi_1,\varphi_2,\varphi_{\alpha_1},\varphi_{\alpha_2}\}$
defines a splitted, unital and factorizable representation of the
unique cycle in $\Gamma$. Hence, there are a map
$u:\{1,\dotsc,m\}\to\{1,\dotsc,m\}$ and scalars $a_1,\dotsc,a_m\in
K$ that satisfy conditions (1) and (2) in Theorem~\ref{th:cycle2}.
Moreover, the maps in $\boldsymbol{\varphi}$ are given by the
formulae \eqref{r1} and \eqref{r2}.

Note that, for an arrow $\alpha$ that is not a loop, we have
$\varphi_{\alpha}=\Id_{K^m}-\varphi_{t(\alpha)}$, as
$\rrank(\Gamma)=1$. Hence, the data that we need to define a
representation of $\Gamma$ is the following:

\begin{definition}\label{dataD}
The data $\mathcal{D}(u,a)$ is defined by
\begin{enumerate}[(1)]
\item
The maps $u_1,u_2,u_3,\dotsc,u_n:\{1,\dotsc,m\}\to\{1,\dotsc,m\}$,
where $u_1=u_2=u$ and for each $i\geq 3$ the maps $u_i$ are
idempotents.
\item
The elements $a_1,\dotsc,a_m\in K$ obeying conditions (1) and (2) in Theorem~\ref{th:cycle2}.
\end{enumerate}
\end{definition}

\begin{theorem}
The data $\mathcal{D}(u,a)$ defines a splitted, unital and
factorizable representation of $\Gamma$ if, and only if, for every
$\alpha\in \Gamma^1\setminus\{\alpha_1,\alpha_2\}$ and every
$p\in\{1,\dotsc,m\}$ such that $u_{s(\alpha)}(p)\neq p$ and
$u_{t(\alpha)}(p)\neq p$, one of the following conditions holds:
\begin{enumerate}[(1)]
\item
$s(\alpha)=1$ and $a_p=1$.
\item
$s(\alpha)=2$ and $a_p=0$.
\end{enumerate}
\end{theorem}
\begin{proof}
Let $\alpha\in \Gamma^1\setminus \{\alpha_1,\alpha_2\}$, and let $p\in\{1,\dotsc,m\}$. We assume that $u_{s(\alpha)}(p)\neq p$ and that $u_{t(\alpha)}(p)\neq p$ and we want to show that either (1) or (2) hold.

There is an arrow $\beta$ such that $t(\beta)=s(\alpha)$, as in the following figure
\[
\xy
 (0,0)*+{{\circ}}="i",(0,3)*+{i},
 (26,0)*+{{\circ}}="j",(26,3)*+{k},
 (13,0)*+{{\circ}}="k",(13,3)*+{j},
 \ar"i";"k"^{\beta}
 \ar"k";"j"^{\alpha}
\endxy
\]
By Lemma~\ref{le:factorizable}(2), we get
\begin{equation}\label{eq:ker}
 \Ker(\varphi_{s(\alpha)})\Ker(\varphi_{t(\alpha)}) = 0.
\end{equation}
First let us prove that $s(\alpha)\in \{1,2\}$. If $s(\alpha)$ is
not a vertex of the cycle, we may apply Lemma~\ref{le:kernel} to
show that $\Ker\varphi_{s(\alpha)}$ is generated by all $f_r$ with
$u_{s(\alpha)}(r)\neq r$. Similarly, $\Ker\varphi_{t(\alpha)}$ is
generated by all $f_r$ with $u_{t(\alpha)}(r)\neq r$. In particular,
by our assumptions, we deduce that $f_p$ is an element in
$\Ker\varphi_{s(\alpha)}\Ker\varphi_{t(\alpha)}$, which is not
possible, in view of \eqref{eq:ker}.

Let us now take an arrow such that $s(a)=1$. We have to prove that
$a_p=1$. Note that \eqref{eq:ker} can be written as follows:
\[
 \Im\varphi_{\alpha_2}\Ker\varphi_k = 0.
\]
Indeed, as $\rank(j)=1$, we have $\beta=\alpha_2$. Thus
$\varphi_j=\Id_A-\varphi_{\alpha_2}$, and the required relation
follows from the fact that $\varphi_{\alpha_2}$ is an idempotent
$K$--linear map. We know that $\Ker(\varphi_k)$ is spanned by all
$f_q$ with $u_k(f_q)\neq f_q$. On the other hand,
\[
 \varphi_{a_2}=\sum_{q=1}^m\left[(a_q - 1)(f_q^\ast-f_{u(q)}^\ast)\right]f_q.
\]
Thus, \eqref{eq:ker} holds if, and only if,
$\varphi_{\alpha_2}(x)f_q=0$, for any $q$ such that $u_k(f_q)\neq
f_q$. At its turn, this property is equivalent to
\[
 (a_q-1)(f_q^\ast-f_{u(q)}^\ast) = 0
\]
for any $q$ such that $u_k(q)\neq q$. As $u(p)\neq p$ and
$u_k(p)\neq p$, we deduce, in particular,  that $a_p=1$. The case
$s(\alpha)=2$ is handled in a similar way.

Conversely, let $\mathcal{D}(u,a)$ be some data as in
Definition~\ref{dataD}. We assume that $\mathcal{D}(u,a)$ satisfy
conditions (1) and (2) of the statement. We define
$\boldsymbol{\varphi}=\{\varphi_1,\varphi_2,\varphi_{\alpha_1},\varphi_{\alpha_2}\}$
by \eqref{r1} and \eqref{r2}. By Theorem~\ref{th:cycle2},
$\boldsymbol{\varphi}$ defines a a splitted unital and factorizable
representation of the cycle in $\Gamma$. Furthermore, for $i\geq 3$
we take $\varphi_i$ to be the algebra endomorphism of $K^m$ that
corresponds to $u_i$. Finally, for an arrow
$\alpha\in\Gamma^1\setminus\{\alpha_1,\alpha_2\}$ that is not a
loop, we set $\varphi_{\alpha}:=\Id_A-\varphi_{t(\alpha)}$.
Obviously, the maps constructed above define a unital and splitted
representation of $\Gamma$. We still have to check that this is
factorizable. Since $\rrank(\Gamma)=1$, in view of
Lemma~\ref{le:factorizable}, we have to prove that
\[
 \Ker(\varphi_{s(\alpha)})\Ker(\varphi_{t(\alpha)})=0
\]
for every arrow $\alpha\in\Gamma^1\setminus\{\alpha_1,\alpha_2\}$.If
$s(\alpha)\in\{1,2\}$, then this equality follows by conditions (1)
and (2) in the statement of the theorem (see the proof of the other
implication). If $s(\alpha)\notin\{\alpha_1,\alpha_2\}$, then the
required relation follows from Lemma~\ref{le:kernel} and the proof
of Theorem~\ref{te:rank1}, taking into account that either
$u_{s(\alpha)}(p)=p$ or $u_{t(\alpha)}(p)=p$.
\end{proof}


\end{document}